\documentclass[12pt]{amsart}
\usepackage{amsmath,amsxtra,amssymb,amsthm,amsfonts}

\newcommand{\limm}{\lim\nolimits}
\newcommand{\limmsup}{\limsup\nolimits}

\newtheorem{theoA}{Theorem}
\newtheorem{theoB}{Theorem}
\newtheorem{theoC}{Theorem}
\newtheorem{corD}{Corollary}

\newcommand{\kla}{\left ( }
\newcommand{\mer}{\right ) }

\renewcommand{\for}{\begin{eqnarray*}}
\newcommand{\mel}{\end{eqnarray*}}
\def\fr{\begin{align*}}

\newcommand{\kl}{\pl \le \pl}
\newcommand{\gl}{\pl \ge \pl}

\newcommand{\lel}{\pl = \pl}

\newcommand{\nz}{{\mathbb N}}
\newcommand{\nen}{n \in \nz}

\newcommand{\rz}{{\mathbb R}}

\newcommand{\cz}{{\mathbb C}}

\newcommand{\ten}{\otimes}

\newcommand{\p}{\hspace{.05cm}}
\newcommand{\pl}{\hspace{.1cm}}

\newcommand{\om}{\omega}

\newcommand{\al}{\alpha}
\newcommand{\si}{\sigma}

\newcommand{\eps}{\varepsilon}

\renewcommand{\L}{{\mathcal L}}

\newcommand{\R}{{\mathcal R}}

\newcommand{\U}{{\mathcal U}}

\newcommand{\summ}{\sum\nolimits}

\newcommand{\prodd}{\prod\nolimits}

\newtheorem{lemma}{Lemma}[section]
\newtheorem{prop}[lemma]{Proposition}
\newtheorem{theorem}[lemma]{Theorem}
\newtheorem{cor}[lemma]{Corollary}

\newtheorem{rem}[lemma]{Remark}

\newcommand{\re}{\begin{rem}\rm}
  \newcommand{\mar}{\end{rem}}

\newcommand{\qd}{\end{proof}\vspace{0.5ex}}

\newcommand{\prf}{\begin{proof}[\bf Proof:]}

\newcommand{\xspace}{\hbox{\kern-2.5pt}}

\allowdisplaybreaks

\begin{document}

\title[Subspaces of noncommutative $L_p$]
{\null \vskip-1cm Rosenthal's theorem \\ for subspaces of
noncommutative \boldmath $L_p$ \unboldmath}

\author[Marius Junge and Javier Parcet]
{Marius Junge and Javier Parcet}

\begin{abstract}
We show that a reflexive subspace of the predual of a von Neumann
algebra embeds into a noncommutative $L_p$ space for some $p>1$.
This is a noncommutative version of Rosenthal's result for
commutative $L_p$ spaces. Similarly for $1 \le q < 2$, an infinite
dimensional subspace $X$ of a noncommutative $L_q$ space either
contains $\ell_q$ or embeds in $L_p$ for some $q < p < 2$. The
novelty in the noncommutative setting is a double sided change of
density.
\end{abstract}

\null

\vskip-30pt

\null

\maketitle

\null

\vskip-25pt

\null

\vspace{-0.3cm}
\section*{Introduction}

The theory of noncommutative $L_p$ spaces has a long tradition in
Banach space theory and the theory of operator algebras
\cite{GK,Haa,Hil,TJ,Fa} and provides the background for recent
progress in noncommutative analysis and probability
\cite{PX1,JLX,JX}. In the commutative setting, the work of
Kadec-Pelczy\'{n}ski \cite{KP} and Rosenthal \cite{Ros} on
subspaces of $L_p$ are corner stones for the understanding of
general Banach space properties. In this paper we prove the
noncommutative version of Rosenthal's result.

\begin{theoA}[Rosenthal '73] A reflexive subspace of $L_1$ embeds
into $L_p$ for some $p>1$.
\end{theoA}

The problem of generalizing Rosenthal theorem to the noncommutative
setting is open for at least 20 years. This problem has an
interesting history. In his seminal paper \cite{P1} on factorization
properties, Pisier described a new approach to some factorization
results by Maurey obtained from Nikishin's theorem. In this paper
Pisier comes very close to proving the noncommutative version of
Rosenthal's result. Indeed, he shows that a reflexive subspace of a
von Neumann algebra predual embeds into an interpolation space
between an $L_1$ space and certain (unusual) $L_2$ space (see
below). Since then it has been a mystery how to modify the argument
and to obtain a subspace of a noncommutative $L_p$ space.
Noncommutative $L_p$ spaces have been defined by Dixmier, Kunze and
Segal in the semifinite setting (see also Nelson \cite{N}) and by
Haagerup \cite{Haa} in the non-tracial case (see also \cite{Hil} for
Connes' approach). Randrianantoanina \cite{Ran} has an argument in
the semifinite setting which is different from ours and does not
provide a good control of the cons\-tants. In this paper we use
modular theory of operator algebras in conjunction with a
noncommutative version of the Peter Jones theorem due to Pisier
\cite{P2} (related to estimates of Kaftal, Larsen and Weiss
\cite{KLW} for triangular matrices) to solve the problem:

\begin{theoB} \label{main0}
Let $N$ be a von Neumann algebra. A reflexive subspace of $L_1(N)$
embeds into $L_p(N)$ for some $p>1$.
\end{theoB}

\vskip1pt

The new interesting point in our proof is the natural change of
density argument. We show that there exists a positive density $d
\in L_1(N)$ such that $tr(d)=1$ and a mapping $u: X \to L_p(N) $
such that
 \[ x \lel d^{1-\frac{1}{p}} u(x) + u(x) d^{1-\frac{1}{p}}
 \pl.\]
In the $\si$-finite case this completely determines $u$. For
simplicity let us assume that $N$ is finite and $d=\sum_j d_je_j$
has a countable spectrum. Then the map $u$ is given by the
following relation
 \[ u(x)\lel \sum_{i,j}^{\null}
 (d_i^{1-\frac1p}+d_j^{1-\frac1p})^{-1} \pl e_ixe_j \pl .\]
Pisier's approach to this result \cite{P1} is used as a  starting
point in our proof. For subspaces of $L_q(N)$ with $q>1$ we have a
similar result, which extends the most general form of Rosenthal's
theorem \cite[Theorem 8]{Ros} to the noncommutative setting.

\begin{theoC} \label{mainq}
Let $N$ be a  von Neumann algebra and fix $1 \le q<2$. Given a
subspace $X$ of $L_q(N)$ not containing $\ell_q$, there exists a
positive density $d \in L_1(N)$ with $tr(d)=1$ and a map $u: X \to
L_p(N)$ for some index $q<p<2$ such that
 \[ x \lel d^{\frac1q-\frac1p} u(x) + u(x) d^{\frac1q-\frac1p}
 \pl.\]
In particular, the space $X$ embeds isomorphically into $L_p(N)$.
\end{theoC}

This result, which also works for linear maps, is closely related
to  Grothendieck type inequalities by Lust-Piquard, see
\cite{lust-gro} and \cite{LPPX}. One of the main obstacles in our
approach to Theorem \ref{mainq} is that the technique of
noncommutative maximal functions is not well-enough understood for
proving Nikishin type results. Therefore we have to work in the
dual setting. Pisier's arguments for $q > 1$ are genuinely very
different from the case $q =1$ which, by duality, leads to linear
maps on $C^*$-algebras. A common characteristic of Pisier's
factorization results in \cite{P1} is a certain differentiation
argument. This is our motivation for the following new inequality.
Let $2 \le p < \infty$ and $a,x$ be positive elements in $L_p(N)$.
Then we have
\begin{equation}\label{diff00}
\|a+x\|_p^p\kl \|a\|_p^p + p \, 2^{p-1} \max \Big\{
\|a^{p-1}x\|_1,\|x\|_p^p \Big\} \pl .
\end{equation}
In the commutative case the triangle inequality in $L_{p-1}$
provides a similar estimate with $2^{p}-1$ instead of $p2^{p-1}$.
For $2\le p\le 3$ operator convexity of $t\mapsto t^{p-1}$
provides an even better estimate. Combined with ultraproduct
techniques, the differential inequality \eqref{diff00} is a
substitute for some of Pisier's arguments in \cite{P1}.

Another technical difficulty concerns complex interpolation of
intersections. We refer to \cite{JP} for many results in this
direction. For a long time, our hope has been to use free
probability to show that interpolation and intersection commute in
this particular setting. However, at the time of this writing some
aspects of harmonic analysis are yet to be discovered before this
approach might be successful. In the Banach space setting of this
paper, we may use different tools from harmonic analysis. Let us
be more specific. We consider a normal faithful state $\phi(\cdot)
= tr(d \, \cdot)$ on a von Neumann algebra $N$ and Pisier's
symmetric norm
 \[ \|x\|_{\Delta_{2}(\phi)} \lel \big( \phi(xx^*) + \phi(x^*x)
 \big)^{\frac12} \sim \max \Big\{ \|d^{\frac12}x\|_2,\|xd^{\frac12}\|_2
 \Big\}  \pl.\]
We will show that
 \begin{equation} \label{maindif}
 \|x\|_{[N,\Delta_{2}(\phi)]_{\frac{2}{p}}} \kl c(p) \max \Big\{
 \|d^{\frac{1}{p}}x\|_p, \|xd^{\frac{1}{p}}\|_p \Big\}
\end{equation}
holds for all $x \in N$ and $2 \le p < \infty$. We can show that
the orthogonal projection from $L_2(N\oplus N)$ to
$\Delta_2(\phi)$ extends to a bounded operator for other values of
$p$.  This allows us to construct the map $u$ in Theorem
\ref{main0}. 

In combination with the results from \cite{JR}, we obtain some
applications to the theory of subsymmetric sequences. A sequence
$(x_n)$ in a Banach space $X$ is called subsymmetric if there
exists a constant $c$ such that $$\Big\| \summ_n a_n x_n \Big\|_X
\sim_c \Big\| \summ_n a_n x_{k_n} \Big\|_X$$ holds for every
strictly increasing sequence $(k_n)$ and arbitrary coefficients
$(a_n)$. We refer to the work of Aldous \cite{Al} and
Krivine-Maurey \cite{Ma-K} for the fact that commutative $L_p$
spaces are stable. This implies in particular that subsymmetric
sequences are symmetric, i.e. we may replace subsequences $(k_n)$
by arbitrary permutation $(\si(n))$. However, due to a result by
Marcolino Nhany  \cite{MN}, noncommutative $L_p$ spaces are in
general not stable.

\begin{corD} \label{subseq}
If $(x_n)\subset N_*$ is a subsymmetric sequence, then $(x_n)$ is
either symmetric or the space $X\lel {\rm span} \big\{x_n \, | \,
n \ge 1 \big\}$ contains $\ell_1$. In particular, $X$ always
contains a symmetric subspace.
\end{corD}

The paper is organized as follows. In section 1 we prove
\eqref{maindif} and the interpolation results for intersections
based on the Peter Jones theorem. This allows us to prove Theorem
\ref{main0} and Corollary \ref{subseq} in section 2. Inequality
\eqref{diff00} and Theorem \ref{mainq} are proved in the last
section of the paper. We use standard notation from the theory of
operator algebras \cite{Tak,KRI,KRII} and the theory of
noncommutative $L_p$ spaces \cite{Terp} (see also
\cite{terp-int}). The reader is assumed to be familiar with basic
ingredients of modular theory and the definition of Haagerup's
noncommutative $L_p$ spaces, see \cite{JX,PX2} for relevant
definitions. However, the main ideas can be understood by \lq
thinking semifinite\rq.

\numberwithin{equation}{section}

\section{An interpolation result}

In this section we provide the main new interpolation results on
intersections and, in particular, the key inequality
\eqref{maindif} will be obtained. In this paper we will use
Haagerup's definition of noncommutative $L_p$ spaces. Indeed, one
first considers the crossed product $M=N\rtimes_{\si_t^{\phi}}\rz$
with respect to a normal semifinite faithful weight $\phi$ on $N$.
Then $M$ is semifinite and there exists a unique normal semifinite
faithful trace $\tau$ on $M$ such that the dual action $\theta_s:
M \to M$ satisfies $\tau(\theta_s(d))=e^{-s}\tau(d)$ for all $s
\in \mathbb{R}$. Haagerup's $L_p$ space is defined as follows
 \[  L_p(N) \lel \Big\{ x \in L_0(M,\tau) \, \big| \ \theta_s(x) =
 e^{-\frac{s}{p}}x \Big\} \pl ,\]
where $L_0(M,\tau)$ stands for the space of $\tau$-measurable
operators affiliated to $M$. For $p=\infty$ we see that
$L_{\infty}(N) = N$ is the set of $\theta_s$-invariant operators.
Let us note that the polar decomposition $x=u|x|$ for $x \in
L_p(N)$ satisfies $u\in N$ and $|x| \in L_p(N)$. We refer to
\cite{Terp} for more details, see also \cite{JX,PX2,TK-II,TK-III}.
An important feature of the Haagerup $L_1(N)$ space is the
distinguished linear functional $tr:L_1(N)\to \cz$. This linear
map implements the isomorphism between $L_1(N)$ and $N_*$. More
precisely, for every normal functional $\phi\in N_*$ there exists
a unique density $d_{\phi}$ with $\phi(x) = tr(d_{\phi}x)$.
Moreover, given $1 \le p < \infty$, the trace functional $tr$ also
implements the duality between $L_p(N)$ and $L_{p'}(N)$. That is,
$L_p(N)^*$ is exactly the space of linear functionals $\phi(x) =
tr(dx)$ with $d \in L_{p'}(N)$ and $\frac1p+\frac{1}{p'}=1$. The
norm in $L_p(N)$ is given by \[ \|x\|_p \lel tr(|x|^p)^{\frac1p}
\pl .\] We also have H\"{o}lder's inequality $\|xy\|_p\kl
\|x\|_q\|y\|_r$ whenever $\frac1p=\frac1q+\frac1r$.

The drawback of Haagerup's construction is the unfamiliar
situation that for $p\neq q$ we have $L_p(N)\cap L_q(N)=\{0\}$. In
particular, this implies that Haagerup $L_p$ spaces do not form an
interpolation scale. However, in this paper interpolation
techniques are important. We shall assume that the reader is
familiar with the complex interpolation method. Let us briefly
review Kosaki's results \cite{Kos} on interpolation of $L_p$
spaces which are crucial in our paper. Once and for all in the
sequel, let us fix a von Neumann algebra $N$ equipped with a
normal faithful state $\phi$ so that $\phi(x) = tr(dx)$ is given
by a positive density $d\in L_1(N)$. Then we may consider the
injective maps
 \[ \iota_{\eta}: x \in N \mapsto d^{1-\eta} x d^{\eta} \in L_1(N)
 \quad \mbox{for} \quad 0 \le \eta \le 1 \pl. \]
 A little bit of modular
theory is required to show that these maps are indeed injective,
see \cite{Kos,JF}. Thus, for fixed $0 \le \eta \le 1$, $(A_0,A_1)
= (\iota_{\eta}(N),L_1(N))$ is an interpolation couple embedded in
$L_1(N)$ as a topological vector space. To be very precise, we
recall that $\|x\|_{A_0} = \|\iota_{\eta}^{-1}(x)\|_N$ and
$\|x\|_{A_1} = \|x\|_{L_1(N)}$. In the literature, the choices
$\eta = 0, \frac12, 1$ are the most important ones. Kosaki showed
that
 \[ \big[ \iota_{\eta}(N), L_1(N) \big]_{\frac1p} \lel
 d^{\frac{1-\eta}{p'}} L_{p}(N) \hskip1pt d^{\frac{\eta}{p'}} \]
holds isometrically. This means exactly that
 \[ \|x\|_p \lel  \big\|
 d^{\frac{1-\eta}{p'}} x d^{\frac{\eta}{p'}}
 \big\|_{[\iota_\eta(N), L_1(N)]_{\frac1p}} \quad \mbox{for all}
 \quad x \in L_p(N) \pl.\]
If $1 \le q < p \le \infty$ and $\frac1s = \frac1q - \frac1p$, we
may also consider the embedding
 \[  \iota_{p,q,\eta}: x \in L_p(N) \mapsto d^{\frac{1-\eta}{s}} x
 d^{\frac{\eta}{s}} \in L_q(N) \pl.\]
Then, the reiteration theorem for complex interpolation
immediately gives
 \[ [\iota_{p,q,\eta}(L_p(N)),L_q(N)]_{\theta}
 \lel \iota_{r,q,\eta}(L_r(N))\]
for $\frac1r = \frac{1-\theta}{p} + \frac{\theta}{q}$. These
interpolation results from \cite{Kos} will be used freely in this
text.

Our aim in this section is to prove a similar result for a double
sided embedding with respect to a fixed density $d$ of a normal
faithful state $\phi$. For $1\le q < p\le \infty$ we define the
following norms
 \[ \|x\|_{L_{p,q}^r(\phi)} \lel \|d^{\frac1q-\frac1p}x\|_q \quad
 \mbox{and} \quad \|x\|_{L_{p,q}^c(\phi)} \lel
 \|xd^{\frac1q-\frac1p}\|_q \pl .\]
Let us write $L_{p,q}^r(\phi)$ and $L_{p,q}^c(\phi)$ for the
respective closures of $L_p(N)$ with respect to the norms above.
Here $r,c$ are chosen because similar expressions appear for
square function inequalities in noncommutative martingale theory.
They correspond to $\eta=0$ and $\eta=1$ in the context of
Kosaki's embedding. We will work with the intersection
 \[ \Delta_{p,q}(\phi) \lel L_{p,q}^r(\phi) \cap L_{p,q}^c(\phi)
 \pl,  \]
defined as the completion of $L_p(N)$ with respect to the norm
 \[ \|x\|_{\Delta_{p,q}(\phi)} \lel \max \Big\{
 \|x\|_{L_{p,q}^r(\phi)}, \|x\|_{L_{p,q}^c(\phi)} \Big\}
 \pl.  \]
Of course, up to an absolute constant, we may replace the maximum
above by the sum or any other $p$-sum. We might use this
equivalence below. We also have a natural embedding
 \[ j_{p,q}:\Delta_{p,q}(\phi)\to L_q(N)\oplus L_q(N) \pl , \quad
 j_{p,q}(x)\lel (d^{\frac1q-\frac1p}x, xd^{\frac1q-\frac1p})
 \pl .\]
According to H\"older's inequality we have a contractive inclusion
$L_p(N) \subset \Delta_{p,q}(\phi)$ given by the identity map.
Therefore $(L_p(N),\Delta_{p,q}(\phi))$ is an interpolation couple
with dense intersection. When $p=\infty$ we shall write
$\Delta_q(\phi)$ for $\Delta_{\infty,q}(\phi)$. Thus, for $(p,q) =
(\infty,2)$ we find the well-known Hilbert space already mentioned
in the Introduction
\[ \|x\|_{\Delta_{2}(\phi)} \sim \kla \|d^{\frac12}x\|_2^2 +
\|xd^{\frac12}\|_2^2 \mer^{\frac12} \lel \phi(xx^* +
x^*x)^{\frac12}
 \lel \sqrt{2} \pl \phi(|x|_s^2)^{\frac12} \pl .\]
Here we followed Pisier's notation $$|x|_s \lel
\sqrt{\frac{x^*x+xx^*}{2}}.$$ Our main result in this section is
the following.

\begin{samepage}
\begin{theorem} \label{minint}
Let $1\le q < p \le \infty$ and
$\frac{1}{r}=\frac{1-\theta}{p}+\frac{\theta}{q}$ for $0 < \theta
< 1$. Then
\begin{itemize}
\item[a)] We have an isomorphism
 \[ \Delta_{p,r}(\phi)\lel
  [L_p(N),\Delta_{p,q}(\phi)]_{\theta}\pl .\]
\item[b)] We may construct a bounded linear map
$$\mathcal{Q}_r: L_r(N\oplus N) \to \Delta_{p,r}(\phi)$$ such that
 \[ \mathcal{Q}_r(d^{\frac1r-\frac1p} x, x d^{\frac1r-\frac1p}) \lel
 x \quad \mbox{for all} \quad x \in L_p(N) \pl. \] In particular,
$j_{p,r}\mathcal{Q}_r$ is a projection from $L_r(N \oplus N)$ onto
$j_{p,r}(\Delta_{p,r}(\phi))$.
\end{itemize}
The relevant constants can be estimated as functions of $p,q,r$ in
both cases.
\end{theorem}
\end{samepage}


\noindent We refer the reader to the end of this section for a
more general form of Theorem \ref{minint}.

\begin{rem} \label{explicit}
{\rm The isomorphism in a)  satisfies
 \[ \Big( \|d^{\frac1r-\frac1p}x\|_r^r +
 \|xd^{\frac1r-\frac1p}\|_r^r
 \Big)^{\frac1r} \sim \|x\|_{[L_p(N),\Delta_{p,q}(\phi)]_{\theta}} \quad
 \mbox{for all} \quad x \in L_p(N) \pl . \]}
\end{rem}

\begin{rem} \label{propertyQ}
{\rm As we shall justify below, the adjoint mapping
$\mathcal{Q}_r^*: \Delta_{p,r}(\phi)^*\to L_{r'}(N\oplus N)$ has
the form $\mathcal{Q}_r^*(\xi)=(u(\xi),u(\xi))$ for some bounded
linear map $u: \Delta_{p,r}(\phi)^*\to L_{r'}(N)$. Equivalently,
we have $\mathcal{Q}_r(y,-y)=0$ for all $y$.}
\end{rem}

It is not very convenient to prove the result for an arbitrary
density $d$. We will apply a well-known construction of Haagerup
and reduce the problem to the case where $N$ is a finite von
Neumann algebra and $d,d^{-1}$ are bounded. Moreover, by
elementary functional calculus, we may then assume that
\begin{equation}\label{simfun}
d \lel \sum_{k=1}^n d_k \hskip1pt e_k
\end{equation}
where the $e_k$ are disjoint projections with $\sum_k e_k=1$ and
$d_k$ are strictly positive numbers such that $d_1 \le d_2 \le
\cdots \le d_n$, see below for justifying this simplification.
Therefore, we will assume in what follows (unless stated
otherwise) that $N$ is finite and that $d$ satisfies
\eqref{simfun}. We note nevertheless that Theorem \ref{minint} is
formulated for Haagerup $L_p$ spaces and hence valid for arbitrary
states. For the moment, we work with a finite von Neumann algebra
and thus we can work with the usual definition of noncommutative
$L_p$ spaces. In particular, all $L_p$ spaces are contained in the
space of $\tau$-measurable operators, see \cite{N} for further
definitions.

In order to sketch our strategy for the proof of part a) in
Theorem \ref{minint}, we need to introduce a more convenient
terminology which will be instrumental in the sequel. Let $d$ be a
density in $L_1(N)$ satisfying \eqref{simfun} and let us write
$L_0(N)$ for the space of $\tau$-measurable operators affiliated
to $N$. Then, given $\al \in \rz$ and $1 \le q \le \infty$, we
define the spaces
\begin{eqnarray*}
L_q^r(N,d^\al) & = & \Big\{x \in L_0(N) \, \big| \ d^\al x \in
L_q(N) \Big\}, \\ L_q^c(N,d^\al) & = & \Big\{x \in L_0(N) \, \big|
\ x d^\al \in L_q(N) \Big\},
\end{eqnarray*}
equipped with the following norms $$\|x\|_{L_q^r(N,d^\al)} =
\|d^\al x\|_q \quad \mbox{and} \quad \|x\|_{L_q^c(N,d^\al)} = \|x
d^\al\|_q \pl .$$ Then, we consider the intersection spaces
$$\Delta_q(N,d^\al) = L_q^r(N,d^\al) \cap L_q^c(N,d^\al)$$ so
that we can recover $\Delta_{p,q}(\phi)$ with $d$ as in
\eqref{simfun} as follows
$$\Delta_{p,q}(\phi) = \Delta_q(N,d^{\frac1q - \frac1p}) \pl .$$
The isomorphism in Theorem \ref{minint} a) is equivalent to
 \begin{equation}
 \label{strat} \Delta_{q_{\theta}}(N,d^{\al_{\theta}})
  \lel [\Delta_{q_0}(N,d^{\al_0}), \Delta_{q_1}(N,d^{\al_1}) ]_{\theta}
 \end{equation}
where $(q_0,q_1,q_\theta) = (p,q,r)$ and $(\al_0,\al_1,\al_\theta)
= (0,1/q - 1/p,\theta \alpha_1)$. As usual we may and will
understand intersections as the diagonal subspaces of the
corresponding direct sum spaces. That is, we have
$$\Delta_{q_j}(N,d^{\al_j}) \subset L_{q_j}^r(N,d^{\al_j})
\oplus L_{q_j}^c(N,d^{\al_j}) \pl, \quad \mbox{for} \quad j=0,1
\pl.$$ By Kosaki's theorem, the components in the direct sum
interpolate isometrically. The easiest way to show that these
intersections commute with interpolation is to show that there is
one projection acting on both spaces $L_{q_j}^r(N,d^{\al_j})
\oplus L_{q_j}^c(N,d^{\al_j})$ for $j=0,1$ which projects onto the
intersection $\Delta_{q_j}(N,d^{\al_j})$. The projection will not
be constructed on $L_{q_j}^r (N,d^{\al_j}) \oplus L_{q_j}^c(N,
d^{\al_j})$ but on spaces of upper and lower triangular elements.

The core of our argument relies on Schur multipliers. This will be
made possible by the canonical embedding $\pi: N \to M_n(N)$ given
by
  \[ \pi(x) \lel \sum_{i,j=1}^n e_{ij} \ten e_i x e_j  \pl .\]
Let us write $\tau$ to denote the trace functional on $N$. This
will allow us to distinguish it from the standard trace $tr$ on
$M_n$. Note that $\pi$ is a normal (not unital) $*$-homomorphism
and we have $(tr \ten \tau) \circ \pi = \tau$. Moreover, the
mapping $\mathsf{E}: M_n(N) \to \pi(N)$ given by $\mathsf{E}(z) =
\pi(1) \hskip1pt z \hskip1pt \pi(1)$ defines a normal conditional
expectation. For the $L_p$-version of the map $\pi$, we first
introduce the normal faithful state $\psi(x)=\sum_k d_k
\tau(x_{kk})$ on $M_n(N)$ with associated density
 \[   \delta = \big( \sum_{k=1}^n d _k e_{kk} \big) \ten 1 \pl
 .\]
With this choice, the map $\pi_p: L_p(N) \to L_p(M_n(N))$
 \begin{equation} \label{ppp}
 \pi_p(d^{\frac{1-\eta}{p}} x d^{\frac{\eta}{p}}) \lel
 \delta^{\frac{1-\eta}{p}} \pi(x) \delta^{\frac{\eta}{p}}
 \end{equation}
becomes an isometric embedding and $\mathsf{E}: L_p(M_n(N)) \to
\pi_p(L_p(N))$ still defines a positive contraction, see \cite{JX}
for further details. Note that $\pi_p^*$ takes $\big( x_{ij} \big)
\in L_{p'}(M_n(N))$ to $\sum_{ij} e_i x_{ij} e_j \in L_{p'}(N)$,
so that $\pi_{p'}^* \pi_p = id_{L_p(N)}$ and $\pi_p \pi_{p'}^* =
\mathsf{E}$. Our main tool are the spaces of lower and upper
triangular matrices in $M_n(N)$ defined as follows
\begin{eqnarray*}
UT_p & = & \Big\{ \big( x_{ij} \big) \in L_p(M_n(N)) \, \big| \
x_{ij}=0 \ \mbox{for $i>j$} \Big\} \pl , \\ LT_p & = & \Big\{
\big( x_{ij} \big) \in L_p(M_n(N)) \, \big| \ x_{ij}=0 \ \mbox{for
$i \le j$} \Big\} \pl.
\end{eqnarray*}
We shall use the fact that $UT_p$ and $LT_p$ are interpolation
scales. This result was proved by Pisier in \cite{P2,P3} and
provides a noncommutative analogue of the Peter Jones theorem on
interpolation of Hardy spaces. We will use the version given in
\cite{PX2}.

\begin{theorem}[Pisier/Xu] \label{upper}
If $1\le p,q \le \infty$ and
$\frac{1}{r}=\frac{1-\theta}{p}+\frac{\theta}{q}$ \[ UT_r \lel
[UT_p,UT_q]_{\theta} \quad \mbox{and} \quad LT_r \lel
[LT_p,LT_q]_{\theta} \pl \] hold with equivalent norms. The
constants are uniformly bounded in $n$.
\end{theorem}

Let us note that for $1<q,p<\infty$ this result follows
immediately from the well-known fact that $UT_p$ and $LT_p$ are
complemented subspaces of $L_p(M_n(N))$. Indeed, the triangular
projection $\mathbf{T}(x_{ij} \ten e_{ij}) = \delta_{i \le j}
\hskip1pt (x_{ij} \ten e_{ij})$ defines a bounded operator on
$L_p(M_n(N))$ with norm controlled by $c\max\{p,p'\}$. Using
$\mathbf{T}$ and $1-\mathbf{T}$ for $p$ and $q$, the interpolation
result follows immediately. The whole point of Pisier's argument
is to extend this result to the non-trivial borderline cases $q=1$
and $p=\infty$.

In our result we are interested in subspaces of $L_p(N)$ which
have upper or lower diagonal form. Moreover, we have to take
different powers of the density $d$ into account. This leads to
consider the following four norms
 \begin{align*}
 \|x\|_{UT_q^{r}(N,d^{\al})} \lel \Big\| \sum_{i\le j}
 d^{\al}_ie_ixe_j \Big\|_{q} &\, ,&
 \|x\|_{UT_q^{c}(N,d^{\al})} \lel \Big\| \sum_{i\le j}
 e_ixe_jd_j^{\al} \Big\|_{q} \pl , \\
 \|x\|_{LT_q^{r}(N,d^{\al})} \lel \Big\| \sum_{i > j}
 d_i^{\al}e_ixe_j \Big\|_{q} &\, ,&
 \|x\|_{LT_q^{c}(N,d^{\al})} \lel \Big\| \sum_{i> j}
 e_ixe_jd_j^{\al} \Big\|_{q} \pl .
 \end{align*}
Then we define the associated spaces
 \begin{eqnarray*}
 UT_q^{r}(N,d^{\al}) & = & \Big\{ x \in L_0(N) \
 \big| \ e_i x e_j = 0 \ \ \mbox{for} \ i > j \, , \
 \|x\|_{UT_q^{r}(N,d^{\al})} < \infty \Big\} \pl , \\
 LT_q^{r}(N,d^{\al}) & = & \Big\{ x \in L_0(N) \
 \big| \ e_i x e_j = 0 \ \ \mbox{for} \ i \le j \, , \
 \|x\|_{LT_q^{r}(N,d^{\al})} < \infty \Big\} \pl , \\
 UT_q^{c}(N,d^{\al}) & = & \Big\{ x \in L_0(N) \
 \big| \ e_i x e_j = 0 \ \ \mbox{for} \ i > j \, , \
 \|x\|_{UT_q^{c}(N,d^{\al})} < \infty \Big\} \pl , \\
 LT_q^{c}(N,d^{\al}) & = & \Big\{ x \in L_0(N) \
 \big| \ e_i x e_j = 0 \ \ \mbox{for} \ i \le j \, , \
 \|x\|_{LT_q^{c}(N,d^{\al})} < \infty \Big\} \pl .
 \end{eqnarray*}
We shall also need to use the spaces
 \begin{eqnarray*}
 UT_q^{r}(M_n(N), \delta^{\al}) & = & \Big\{ \big(x_{ij} \big)
 \in L_0(M_n(N)) \ \big| \ x_{ij} = 0 \ \ \mbox{for} \ i > j
 \, , \ \|\delta^\al (x_{ij})\|_q < \infty \Big\} \pl
 , \hskip3pt \null \\ LT_q^{r}(M_n(N), \delta^{\al}) & = &
 \Big\{ \big( x_{ij} \big) \in  L_0(M_n(N)) \ \big| \ x_{ij} = 0
 \ \ \mbox{for} \ i \le j \, , \ \|\delta^\al (x_{ij})\|_q <
 \infty \Big\} \pl , \\ UT_q^{c}(M_n(N), \delta^{\al}) & = &
 \Big\{ \big( x_{ij} \big) \in L_0(M_n(N)) \ \big| \ x_{ij} = 0 \
 \ \mbox{for} \ i > j \, , \ \|(x_{ij}) \delta^\al\|_q < \infty
 \Big\} \pl , \\ LT_q^{c}(M_n(N), \delta^{\al}) & = & \Big\{
 \big( x_{ij} \big) \in L_0(M_n(N)) \ \big| \ x_{ij} = 0 \ \
 \mbox{for} \ i \le j \, , \ \|(x_{ij}) \delta^\al\|_q < \infty
 \Big\} \pl .
 \end{eqnarray*}
Let us observe that, if $e_i x e_j = 0$ for $i
> j$, we have for $\alpha = \frac1q - \frac1p$
$$\pi_q \Big( \sum_{i \le j} d_i^\alpha e_i x e_j \Big) = \pi_q(d^\alpha
x) = \delta^\alpha \pi_p(x) \pl .$$ In particular, it is easily
seen that
\begin{equation} \label{projUT}
\mathsf{E}: UT_q^r(M_n(N), \delta^\al) \to \pi_p \big( UT_q^r(N,
d^\al) \big)
\end{equation}
is still a contractive projection. This property (which extends
automatically to the three other spaces considered above) will be
instrumental in the following result, where we combine Kosaki's
embedding with interpolation of triangular matrices.

\begin{lemma} \label{intpis0}
If $1 \le q_0, q_1 \le \infty$ and $\alpha_0, \alpha_1 \in
\mathbb{R}$, let us take $1/q_\theta = (1-\theta)/q_0 +
\theta/q_1$ and $\alpha_\theta = (1-\theta) \al_0 + \theta \al_1$.
Then, the following isomorphisms hold with relevant constants
depending only on $q_0, q_1$ and $\theta$
\begin{eqnarray*}
\big[ UT_{q_0}^r (N,d^{\al_0}), UT_{q_1}^r (N,d^{\al_1})
\big]_{\theta} & = & UT_{q_{\theta}}^r (N,d^{\al_\theta}) \pl , \\
\big[ LT_{q_0}^r (N,d^{\al_0}) \hskip1pt , \hskip1pt LT_{q_1}^r
(N,d^{\al_1})
\big]_{\theta} & = & LT_{q_{\theta}}^r (N,d^{\al_\theta}) \pl , \\
\big[ UT_{q_0}^c (N,d^{\al_0}), UT_{q_1}^c (N,d^{\al_1})
\big]_{\theta} & = & UT_{q_{\theta}}^c (N,d^{\al_\theta}) \pl , \\
\big[ LT_{q_0}^c (N,d^{\al_0}) \hskip1pt , \hskip1pt LT_{q_1}^c
(N,d^{\al_1}) \big]_{\theta} & = & LT_{q_{\theta}}^c
(N,d^{\al_\theta}) \pl .
\end{eqnarray*}
\end{lemma}

\begin{proof}
Since the proof of the four isomorphisms is identical, we only
consider the first one. According to the boundedness of
\eqref{projUT}, it suffices to prove the analogous isomorphism on
the amplified algebra $M_n(N)$
 \begin{equation}\label{large}
 \big[ UT_{q_0}^r(M_n(N), \delta^{\al_0}), UT_{q_1}^r
 (M_n(N),\delta^{\al_1}) \big]_{\theta}
 \lel UT_{q_{\theta}}^r (M_n(N),\delta^{\al_\theta})  \pl . \end{equation}
Indeed, Kosaki's interpolation theorem tells us that
 \[ \null \hskip2pt
 \big[ L_{q_0}^r(M_n(N),\delta^{\al_0}),L_{q_1}^r(M_n(N),
 \delta^{\al_1}) \big]_{\theta}
 \lel L_{q_{\theta}}^r(M_n(N),\delta^{\al_\theta})  \]
holds isometrically. Thus, by our special choice of $\delta$, we
obtain a contractive inclusion
 \[ \big[ UT_{q_0}^r(M_n(N),\delta^{\al_0}), UT_{q_1}^r (M_n(N),
 \delta^{\al_1}) \big]_{\theta} \, \subset \,
 UT_{q_{\theta}}^r(M_n(N),\delta^{\al_\theta}) \pl. \]
For the converse, we assume that $x \in UT_{q_\theta}^r(M_n(N),
\delta^{\alpha_\theta})$ has norm less than $1$. That is, $x \in
L_0(M_n(N))$ is an upper triangular matrix such that
$\|\delta^{\al_\theta}x\|_{q_\theta}< 1$. Let $\mathcal{S}$ stand
for the strip $\mathcal{S} = \big\{ z \in \mathbb{C} \, \big| \ 0
\le \mbox{Re}(z) \le 1 \big\}$ and denote by
$(\partial_0,\partial_1)$ the left and right sides of its
boundary. According to Theorem \ref{upper}, we may find an
analytic function $$f: \mathcal{S} \to UT_{q_0} + UT_{q_1}$$ such
that $f(\theta)=\delta^{\al_\theta} x$ and
 \[ \max \Big\{ \sup_{z \in
 \partial_0} \|f(z)\|_{UT_{q_0}}, \, \sup_{z \in \partial_1}
 \|f(z)\|_{UT_{q_1}} \Big\} \le c(q_{\theta})\]
holds for some universal constant $c(q_{\theta})$. Then we define
$g(z) \lel \delta^{-(1-z)\al_0-z\al_1} f(z)$. Note that $g$ is
analytic and that $g(z)$ is still an upper triangular matrix for
any $z \in \mathcal{S}$. For $z \in
\partial_0$ we find
 \[ \|g(z)\|_{UT_{q_0}^r(M_n(N),\delta^{\al_0})} \lel \|f(z)
  \|_{UT_{q_0}} \kl c(q_{\theta})
 \pl .\]
Similarly, if $z \in \partial_1$ we have the estimate
 \[ \|g(z)\|_{UT_{q_1}^r(M_n(N),\delta^{\al_1})} \lel
  \|f(z)\|_{UT_{q_1}} \kl c(q_{\theta}) \pl .\]
Clearly we have $g(\theta)\lel x$ and \eqref{large} follows from
the three lines lemma. \qd

The next lemma is a very well-known classical result. We have
decided to include the proof for the convenience of the reader.
The easy argument that we use here is due to Burak Erdogan.

\begin{lemma}\label{pd} Let $f: \mathbb{R} \to \mathbb{R}$ be an
even integrable function whose restriction to $\mathbb{R}_+$ is
non-increasing and convex. Assume that $f$ is differentiable
almost everywhere and $f'$ is integrable. Then $f$ is positive
definite, i.e. its Fourier transform is positive.
\end{lemma}

\begin{proof} If $\xi \in \mathbb{R}_+$, we have
 \begin{eqnarray*}
 \widehat{f}(\xi) & = & \int_{\mathbb{R}} f(x) \hskip1pt
 e^{-i x \xi} dx \lel 2 \int_{\mathbb{R}_+} f(x) \cos(x\xi)
 \hskip1pt dx \\ & = & - \frac{2}{\xi} \int_{\mathbb{R}_+} f'(x)
 \sin(x \xi) \hskip1pt dx \lel - \frac{2}{\xi^2} \int_{\mathbb{R}_+}
 f'(\frac{x}{\xi}) \sin(x) \hskip1pt dx
 \pl .
 \end{eqnarray*}
Here we used the fact that $f$ is even, integration by parts and
substitution. The function $g(x)=-f'(\frac{x}{\xi})$ is positive,
non-increasing and integrable on $\mathbb{R}_+$. In particular, we
deduce that
 \[ \gamma_k = \int_0^{2 \pi} g(x + 2 \pi k) \sin(x)
 \hskip1pt dx \gl 0 \] for all integer $k \ge 0$ and therefore
 \[ \widehat{f}(\xi) \lel \frac{2}{\xi^2} \sum_{k \ge 0}
 \gamma_k \gl 0 \quad \mbox{for all} \quad x \in \mathbb{R}_+
 \pl. \]
By symmetry, $\widehat{f}(\xi) \ge 0$ for all $\xi \neq 0$.
Moreover, since $f$ is positive, we have $$\widehat{f}(0) =
\int_{\mathbb{R}} f(x) \hskip1pt dx \gl 0.$$ This shows that
$\widehat{f}: \mathbb{R} \to \mathbb{R}_+$, so that $f$ is
positive definite and the proof is complete. \qd

\begin{lemma} \label{Schur}
Let $a = \big( \sum_k a_k e_{kk} \big) \ten 1$ be a positive
density on $M_n(N)$ with non-decreasing entries $a_1 \le a_2 \le
\cdots \le a_n$. Let $\mathcal{L}_a(x) = ax$ and $\R_a(x) = xa$ be
the left and right multiplication maps. Then, the norm of the maps
\[ \mathcal{L}_{a^{\eta}} \R_{a^{1-\eta}} (\mathcal{L}_{a} +
 \R_a)^{-1} \quad (0 \le \eta \le 1) \]
on the spaces $UT_p$ and $LT_p$ is bounded by $\frac32$ for all $1
\le p \le \infty$. In particular, given $\al,\beta \in \rz$ and
$d$ a density as in \eqref{simfun}, the norm of the following maps
is also bounded by $\frac32$ on $UT_q^{r/c}(N,d^{\al})$ and
$LT_q^{r/c}(N,d^{\al})$ for all $1 \le q \le \infty$ and all $0
\le \eta \le 1$
 \[ \L_{d^{(1-\eta)\beta}} \R_{d^{\eta\beta}}
 (\L_{d^{\beta}}+ \R_{d^{\beta}})^{-1} \pl .\]
\end{lemma}

\begin{proof}
Let $x\in  UT_p$ be an upper triangular matrix.
Then we observe that
\[ \mathcal{L}_{a} (\mathcal{L}_{a} + \R_a)^{-1}(x_{ij}) \lel
 \Big( \frac{a_i}{a_i+a_j} \, x_{ij} \Big) \lel
 \Big( \frac{\min(a_i,a_j)}{a_i+a_j} \, x_{ij} \Big) \]
because for $i>j$ we have $x_{ij}=0$. Observe that the same
argument shows that on $LT_p$ we have to use $\max(a_i,a_j)$
instead of $\min(a_i,a_j)$. However, we have
$\frac{\max(a_i,a_j)}{a_i+a_j} = 1 -
\frac{\min(a_i,a_j)}{a_i+a_j}$. Therefore, the cases $\eta = 0,1$
follow immediately once we have shown that
 \[ M_a(x_{ij}) \lel \Big( \frac{\min(a_i,a_j)}{a_i+a_j} \, x_{ij}
 \Big) \]
is bounded on $L_p(M_n(N))$ for all $1\le p\le \infty$. If $s,t
\in \mathbb{R}_+$, we have
 \[ \frac{\min(s,t)}{s+t}
 \lel \frac{1}{1+\frac{\max(s,t)}{\min(s,t)}}
 \lel \frac{1}{1+e^{|\log(s) - \log(t)|}} \pl .\]
The Fourier inversion formula for $f(x)=\frac{1}{1+e^{|x|}}$ gives
$$\frac{\min(s,t)}{s+t} \lel
  \frac{1}{1+e^{|\log(s)-\log(t)|}} \lel
  \frac{1}{2 \pi}  \int_\rz \widehat{f}(\xi)
  \hskip1pt e^{i\xi(\log(s)-\log(t))} \hskip1pt d \xi \pl.$$
According to Lemma \ref{pd}, $f$ is positive definite and we
obtain
 \begin{eqnarray*}
  \|M_a(x_{ij})\|_p & = & \Big\| \Big( \frac{1}{2 \pi} \int_\rz
  \widehat{f}(\xi) \hskip1pt e^{i \xi(\log(a_i) - \log(a_j))} x_{ij}
  \hskip1pt d\xi \Big) \Big\|_p \\ & \le & \frac{1}{2 \pi}
  \int_{\mathbb{R}} \widehat{f}(\xi) \Big\| \big( e^{i \xi \log(a_i)}
  x_{ij} e^{-i \xi \log(a_j)} \big) \Big\|_p \hskip1pt d\xi \\
  & \le & \frac{1}{2\pi} \int_{\mathbb{R}} \widehat{f}(\xi) \hskip1pt
  d\xi \, \|(x_{ij})\|_p \lel \frac{\|(x_{ij})\|_p}{1+e^0} =
  \frac12 \, \|(x_{ij})\|_p \pl .
  \end{eqnarray*}
Thus, $M_a$ is bounded on $L_p(M_n(N))$ with norm $\frac12$ and
the same holds for $\L_{a} (\L_{a} + \R_a)^{-1}$ on the space
$UT_p$. Moreover, the same arguments show that $\R_{a} (\L_{a} +
\R_a)^{-1}$ is bounded on $LT_p$ with norm $\frac12$. On the other
hand, recalling one more time that $\frac{\max(a_i,a_j)}{a_i +
a_j} + \frac{\min(a_i,a_j)}{a_i + a_j} = 1$, we deduce that
$\L_{a} (\L_{a} + \R_a)^{-1}$ on $LT_p$ and $\R_{a} (\L_{a} +
\R_a)^{-1}$ on $UT_p$ are respectively bounded by $1 + \frac12$.
It remains to prove the case $0 < \eta < 1$. Let us consider $x
\in UT_p$ and define the complex function $f(z) \lel \L_{a^{1-z}}
\R_{a^z} (\L_a + \R_a)^{-1}(x)$. Then it is easily seen that
 \[ \max \Big\{ \sup_{z \in
 \partial_0} \|f(z)\|_p, \, \sup_{z \in \partial_1} \|f(z)\|_p \Big\}
 \le \frac32 \, \|x\|_p \pl . \]
Thus, we find that $\|f(\eta)\|_p \le \frac32 \|x\|_p$. The
argument for $LT_p$ is similar. Let us now prove the second
assertion. Since the left and right multiplication maps
$\mathcal{L}$ and $\mathcal{R}$ clearly commute with $d^\al$, it
is no restriction to assume that $\al = 0$ and $q=p$. On the other
hand, taking $$z_{ij} = d_i^{(1-\eta)\beta} \frac{e_i x
e_j}{d_i^\beta + d_j^\beta} \, d_j^{\eta \beta} \pl ,$$ we clearly
have $$\L_{\delta^{(1-\eta)\beta}}\R_{\delta^{\eta\beta}}
(\L_{\delta^{\beta}}+ \R_{\delta^{\beta}})^{-1} \pi_p(x) =
\summ_{ij} e_{ij} \ten z_{ij} = \pi_p \L_{d^{(1-\eta)\beta}}
\R_{d^{\eta\beta}} (\L_{d^{\beta}}+ \R_{d^{\beta}})^{-1}(x) \pl
.$$ Therefore, the first assertion implies the second assertion
and we are done. \qd

\noindent In the following we use the notations
 \begin{eqnarray*}
  UT_q(N) & = & UT_q^r(N,d^{0}) \lel UT_q^c(N,d^{0}) \pl , \\
  LT_q(N) & = & \hskip1pt LT_q^r(N,d^{0}) \lel \hskip1pt
  LT_q^c(N,d^{0}) \pl , \\ \Delta_q^{UT}(N,d^{\al}) & = &
  UT_q^r(N,d^{\al}) \cap  UT_q^c(N,d^{\al}) \pl , \\
  \Delta_q^{LT}(N,d^{\al}) & = &
  LT_q^r(N,d^{\al}) \hskip1pt \cap \hskip1pt LT_q^c(N,d^{\al})
  \pl ,
  \end{eqnarray*}
for spaces of upper and lower triangular elements.

\begin{lemma}\label{proj} Let $1 \le q_0, q_1 \le \infty$,
$\al\in \rz$ and $\al_\theta=\theta\al$. Then the map
 \[ \Lambda: UT_{q_{\theta}}^r(N,d^{\al_{\theta}})
 \oplus UT_{q_{\theta}}^c(N,d^{\al_{\theta}}) \to
   \big[ UT_{q_0}(N),
  \Delta_{q_1}^{UT}(N,d^{\al}) \big]_{\theta} \]
defined by
 \[  \Lambda(y,z) \lel
 (\L_{d^{\al_\theta}}+\R_{d^{\al_\theta}})^{-1}(d^{\al_\theta} y
 + z d^{\al_\theta})\]
satisfies $\|\Lambda\|\kl c(q_{\theta})$. The same holds for the
space of lower triangular matrices.
\end{lemma}

\begin{proof} According to Lemma \ref{intpis0}, we know that
\begin{eqnarray*}
UT_{q_{\theta}}^r(N,d^{\al_{\theta}}) & = & \big[ UT_{q_0}(N),
UT_{q_1}^r(N,d^{\al}) \big]_{\theta} \pl ,
\\ UT_{q_{\theta}}^c(N,d^{\al_{\theta}}) & = &
\big[ UT_{q_0}(N), UT_{q_1}^c(N,d^{\al}) \big]_{\theta} \pl ,
\end{eqnarray*}
holds up to a constant $c'(q_{\theta})$. Obviously, we have
$\Lambda(x,x)=x$. Therefore, it suffices to show that $\Lambda$ is
bounded on $UT_{q_0}(N) \oplus UT_{q_0}(N)$ and on
$UT_{q_1}^r(N,d^{\al}) \oplus UT_{q_1}^c(N,d^{\al})$. Indeed, we
deduce from Lemma \ref{Schur} that
 \begin{eqnarray*}
 \|\Lambda(y,z)\|_{UT_{q_0}(N)}
 & = & \big\| \L_{d^{\al}}(\L_{d^{\al}}+\R_{d^{\al}})^{-1}(y)+
 \R_{d^{\al}}(\L_{d^{\al}} + \R_{d^{\al}})^{-1}(z) \big\|_{q_0} \\
 & \le & \frac32 \, \|y\|_{q_0} + \frac32 \, \|z\|_{q_0} \kl
 3 \pl \|(y,z)\|_{UT_{q_0}(N) \oplus
 UT_{q_0}(N)} \pl .
 \end{eqnarray*}
On the other hand, we have
 \begin{eqnarray*}
 \|\Lambda(y,z)\|_{UT_{q_1}^r(N,d^{\al})}
 & = & \big\| \L_{d^{\al}} (\L_{d^{\al}} + \R_{d^{\al}})^{-1}(d^{\al}y +
 zd^{\al}) \big\|_{q_1} \kl \frac32 \, \big\|d^{\al}y + zd^{\al} \big\|_{q_1}
 \\ [4pt] & \le & 3 \pl
 \max \Big\{ \|d^{\al}y\|_{q_1}, \|zd^{\al}\|_{q_1} \Big\} \lel
 3 \pl \|(y,z)\|_{UT_{q_1}^r(N,d^{\al}) \oplus
 UT_{q_1}^c(N,d^{\al})} \pl .
 \end{eqnarray*}
The estimate for $UT_{q_1}^c(N,d^{\al})$ uses
$\R_{d^{\al}}(\L_{d^{\al}}+\R_{d^{\al}})^{-1}$ instead. On the
other hand, the proof for lower triangular matrices is verbatim
the same. The proof is complete. \qd

The next result is well-known. It can be proved using the fact
that $L_p(N)$ are UMD spaces  (see \cite{BerkG} and \cite{Bou}) or
applying the boundedness for the noncommutative Hilbert transform
in chapter 8 of \cite{PX2}, see also the earlier results in
\cite{GK,KPT}.

\begin{lemma} \label{tri}
Let $(e_i)$ be a family of disjoint projections in a von Neumann
algebra $N$ and let us consider the triangular projection
$\mathbf{T}_e(x) = \sum_{i \le j} e_i x e_j$. Then, the mapping
$\mathbf{T}_e$ is bounded on $L_p(N)$ for $1<p<\infty$.
\end{lemma}

\begin{proof} It is well-known that the triangular projection
$\mathbf{T}(a) = \sum_{i\le j} a_{ij} e_{ij}$ is completely
bounded on $S_p$, see the references above.  Then, the bounded map
$\pi_{p'}^* \mathbf{T} \pi_p$ yields the modified triangular
projection $\mathbf{T}_e$ used in the assertion. \qd

\begin{proof}[{\bf Step 1 of the proof.}] We will
prove Theorem \ref{minint} assuming (1.1). For the first assertion
a), we observe from Kosaki's interpolation that the inclusion
$[L_p(N), \Delta_{p,q}(\phi)]_\theta \subset \Delta_{p,r}(\phi)$
is trivially contractive. For the converse we use the
$\Delta_q(N,d^\al)$ terminology. In other words we have to prove
that $$\Delta_r(N,d^{\theta/s}) \subset \big[ L_p(N), \Delta_q(N,
d^{1/s}) \big]_\theta \quad \mbox{with} \quad \frac1s = \frac1q -
\frac1p \pl .$$ On the other hand, the inclusions $UT_p(N) \subset
L_p(N)$ and $\Delta_q^{UT}(N,d^{1/s}) \subset \Delta_q(N,d^{1/s})$
are contractive and the same happens for the spaces of lower
triangular matrices. Therefore, considering the decomposition $x =
\mathbf{T}_e(x) + x - \mathbf{T}_e(x)$ for $x \in
\Delta_r(N,d^{\theta/s})$, it suffices to show that
\begin{eqnarray*}
\Delta_r^{UT}(N,d^{\theta/s}) & \subset & \big[ UT_p(N),
\Delta_q^{UT}(N, d^{1/s}) \big]_\theta \pl , \\
\Delta_r^{LT}(N,d^{\theta/s}) & \subset & \big[ LT_p(N) \hskip1pt,
\hskip1pt \Delta_q^{LT}(N, d^{1/s}) \big]_\theta \pl .
\end{eqnarray*}
Note that $1 < r < \infty$ because $0<\theta<1$. According to
Lemma \ref{tri}, this implies that $\mathbf{T}_e(x)$ belongs to
$\Delta_r^{UT}(N,d^{\theta/s})$ and $x - \mathbf{T}_e(x) \in
\Delta_r^{LT}(N,d^{\theta/s})$. Hence, applying Lemma \ref{proj}
we deduce that
 \begin{eqnarray*}
 \|\mathbf{T}_e(x)\|_{[UT_{p}(N),\Delta_q^{UT}(N,d^{1/s})]_{\theta}}
 & = & \|\Lambda(\mathbf{T}_e(x),\mathbf{T}_e(x))\|_{[UT_{p}(N),
 \Delta_q^{UT}(N,d^{1/s})]_{\theta}} \\
 & \le &  c(r) \|\mathbf{T}_e(x)\|_{\Delta_r^{UT}(N,d^{\theta/s})} \kl
 c(r) d(r) \|x\|_{\Delta_r(N,d^{\theta/s})} \pl .
 \end{eqnarray*}
The same argument with respect to lower triangular matrices gives
 \[ \|x-\mathbf{T}_e(x)\|_{[LT_{p}(N),\Delta_q^{LT}
 (N,d^{1/s})]_{\theta}} \kl c(r) d(r) \|x\|_{\Delta_r(N,d^{\theta/s})}
 \pl . \]
For the proof of part b) we construct
 \[ \mathcal{Q}_r: L_r(N\oplus N) \to \Delta_{p,r}(\phi) \]
as follows
 \begin{equation}
 \label{expff}
  \mathcal{Q}_r(y,z)
  \lel (\L_{d^{\frac1r-\frac1p}}+\R_{d^{\frac1r-\frac1p}})^{-1}
  (y+z) \pl
  \end{equation}
for $y,z\in L_r(N)$. Clearly, we have
$\mathcal{Q}_r(d^{\frac1r-\frac1p} x, x d^{\frac1r-\frac1p}) = x$
for all $x \in L_p(N)$. For the norm estimate, we use again the
fact that the triangular map $\mathbf{T}_e$ is bounded. For $y,z
\in L_r(N)$ we deduce from Lemma \ref{Schur} that
 \begin{eqnarray*}
 \big\| \mathcal{Q}_r \big( \mathbf{T}_e(y), \mathbf{T}_e(z) \big)
 \big\|_{UT_r^r(N,d^{\frac1r-\frac1p})}
 & = & \big\| \L_{d^{\frac1r-\frac1p}}
 (\L_{d^{\frac1r-\frac1p}}+\R_{d^{\frac1r-\frac1p}})^{-1}
 \big( \mathbf{T}_e(y) + \mathbf{T}_e(z) \big) \big\|_r \\
 & \le & 3 \, \big\| \big( \mathbf{T}_e(y), \mathbf{T}_e(z) \big)
 \big\|_{L_r(N \oplus N)} \kl 3 \, d(r) \|(y,z)\|_{L_r(N \oplus N)}
 \pl ,
 \end{eqnarray*}
where $d(r)$ stands for the norm of the triangular projection on
$L_r(N)$. The same estimate holds for
$UT_r^c(N,d^{\frac1r-\frac1p})$. We can also repeat the estimate
for $y-\mathbf{T}_e(y)$ and $z-\mathbf{T}_e(z)$ with respect to
the spaces $LT_r^r(N,d^{\frac1r-\frac1p})$ and
$LT_r^c(N,d^{\frac1r-\frac1p})$. This yields the norm estimate
 \[ \big\| \mathcal{Q}_r: L_r(N\oplus N) \to \Delta_{p,r}(\phi)
 \big\| \kl 6 \, d(r)
 \pl . \qedhere \] \end{proof}

\begin{rem}{\rm In our applications we will combine a) and b) and
deduce that $$\|\mathcal{Q}_r:L_r(N\oplus N) \to [L_p(N),
\Delta_{p,q}(\phi)]_{\theta} \|\kl 6 \, c(r) \, d(r)^2 \pl .$$
Here we use the triangular projection twice. In the first version
of this paper we directly constructed a map
$\widehat{\mathcal{Q}}_r: L_r(N\oplus N) \to
[L_p(N),\Delta_{p,q}(\phi)]_{\theta}$ projecting onto the
canonical image of $\Delta_{p,r}(\phi)$. Indeed, we consider upper
triangular elements $y = d^{\frac1r-\frac1p} x_1, z = x_2
d^{\frac1r-\frac1p}$ with $x_1,x_2\in UT_p(N)$. Then the
``canonical'' image in $L_q(N\oplus N)$ is given by
 \[ (d^{\frac1q-\frac1r} y , z d^{\frac1q-\frac1r})
 \lel (d^{\frac1q - \frac1p} x_1, x_2 d^{\frac1q - \frac1p}) \pl .\]
We have seen in Lemma \ref{proj} that
 \[ \|\Lambda(x_1,x_2)\|_{\Delta_q^{UT}(N,d^{\al})}
 \kl 3 \pl \big\| (d^{\frac1q-\frac1r} y , z d^{\frac1q-\frac1r})
 \big\|_{L_q(N\oplus N)} \quad \mbox{for} \quad \alpha = 1/q-1/p
 \pl.\]
The same estimate holds with respect to $LT_p(N)$. By complex
interpolation we deduce
 \begin{eqnarray*}
 \lefteqn{\|\Lambda(x_1,x_2)\|_{[UT_p(N),\Delta_q^{UT}
 (N,d^{\frac1q-\frac1p})]_{\theta}}} \\ & \le & 3 \pl
 \|(x_1,x_2)\|_{[UT_p(N) \oplus
 UT_p(N), UT_{q}^r(N,d^{\frac1q-\frac1p}) \oplus
 UT_{q}^c(N,d^{\frac1q-\frac1p})]_{\theta}}\\
 & \le & 3 \pl c(r) \pl \|(x_1,x_2)\|_{UT_r^r
 (N,d^{\frac1r-\frac1p}) \oplus UT_r^c(N,d^{\frac1r-\frac1p})}
 \lel 3 \pl c(r) \pl \|(y,z)\|_{L_r(N\oplus N)}
 \pl .
 \end{eqnarray*}
For $y,z\in L_r(N)$ we consider $\zeta= d^{\frac1q-\frac1r} y + z
d^{\frac1q-\frac1r}$ and the projection
 \begin{align*}
  \widehat{\mathcal{Q}}_r(y,z)&=
  \Big( d^{\frac1q-\frac1p} \big((\L_{d^{\frac1q-\frac1p}}
  + \R_{d^{\frac1q-\frac1p}})^{-1} (\zeta)\big),
  \big((\L_{d^{\frac1q-\frac1p}}+\R_{d^{\frac1q-\frac1p}})^{-1}
  (\zeta)\big) d^{\frac1q-\frac1p} \Big) \pl.
 \end{align*}
Then we have
 \[ \big\| \widehat{\mathcal{Q}}_r: L_r(N\oplus N) \to [L_p(N),
 \Delta_{p,q}(\phi)]_{\theta} \big\| \kl 6 \pl c(r) d(r) \pl .\]
It is known that $d(r) \kl c \max\{r,r'\}$. However, we have no
explicit control on $c(r)$. It would be interesting to know
whether the singularity for $r \to 1$ is necessary when
interpolating $[N,\Delta_{\infty,1}(\phi)]_{\frac1r}$.}
\end{rem}

\begin{rem}
{\rm In contrast to $\widehat{\mathcal{Q}}_r$, the projection from
part b) satisfies the condition $\mathcal{Q}_r(y,-y)=0$ mentioned
in Remark \ref{propertyQ}. This follows immediately from
\eqref{expff} and is important for our applications below. Let us
reformulate this condition for the dual map. Using
$\mathcal{Q}_r(y,-y)=0$ we see that $\mathcal{Q}_r$ factors
through
 \[ L_r(N) \simeq L_r(N\oplus N)/\{(y,-y) \, | \ y\in
 L_r(N)\}\pl .\]
More explicitly, for $\xi\in \Delta_{p,r}(\phi)^*$ we have
 \begin{eqnarray*}
 \big\langle \mathcal{Q}_r^*(\xi),(y,z) \big\rangle & = &
 \big\langle \xi, \mathcal{Q}_r(y,z) \big\rangle \lel \big\langle
 \xi, (\L_{d^{\frac1r-\frac1p}} + \R_{d^{\frac1r-\frac1p}})^{-1}
 (y+z) \big\rangle \\ & = & \big\langle \xi,
 (\L_{d^{\frac1r-\frac1p}} + \R_{d^{\frac1r-\frac1p}})^{-1}(y)
 \big\rangle + \big\langle \xi,
 (\L_{d^{\frac1r-\frac1p}} + \R_{d^{\frac1r-\frac1p}})^{-1}(z)
 \big\rangle \pl .
 \end{eqnarray*}
This allows us to define the bounded map $u(\xi)$ by $\langle
u(\xi),y\rangle=\frac12 \langle
\mathcal{Q}_r^*(\xi),(y,y)\rangle$. Clearly, we have
$\mathcal{Q}_r^*(\xi)=(u(\xi),u(\xi))$. Assuming \eqref{simfun}
the map $(\L_{d^{\frac1r-\frac1p}} +
\R_{d^{\frac1r-\frac1p}})^{-1}$ is bounded. In the next steps of
our proof this is not necessarily the case, but see Corollary
\ref{15} below.}
\end{rem}

\begin{proof}[{\bf Step 2 of the proof.}]
We now study the case where $N$ is finite and equipped with a
density $d$ such that $c_1 1 \le d \le c_2 1$ for some constants
$0 < c_1 \le c_2 < \infty$, so that $d$ and $d^{-1}$ are bounded.
We claim that for any $\eps>0$ we may find a density $d_\eps$ of
the form \eqref{simfun}, with $\tau(d_\eps) = 1$ and such that
$$(1+\eps)^{-1} d_{\eps} \le d \le (1+\eps) \hskip1pt d_{\eps} \pl
.$$ Indeed, let $\mu$ be the probability measure on the Borel
$\sigma$-algebra over $[c_1,c_2]$ determined by $\mu(\mathsf{E}) =
\tau(1_{\mathsf{E}}(d))$, where $1_{\mathsf{E}}(d)$ denotes the
corresponding spectral projection. This provides isometric
isomorphisms $L_p(\mu) = L_p(A,\tau)$, where $A$ is the (abelian)
von Neumann subalgebra of $N$ given by
$$A = \Big\{ f(d) \, \big| \ f: [c_1,c_2] \to \mathbb{C} \
\mbox{bounded and measurable} \Big\}.$$ In particular, we may
approximate $d$ by $d_\eps$ of the form \eqref{simfun} just by
approximating the function $f(x) = x$ by a suitable simple
function. In particular, we may even assume that $d_{\eps}$
commutes with $d$. Letting $\phi_\eps(x) = tr(d_\eps x)$ be the
state determined by $d_\eps$ and taking $\frac1s = \frac1q -
\frac1p$, it is clear that
$$(1+\eps)^{\frac{-1}{s}} \|x\|_{\Delta_{p,q}(\phi_{\eps})} \kl
\|x\|_{\Delta_{p,q}(\phi)} \kl (1+\eps)^{\frac1s} \hskip1pt
\|x\|_{\Delta_{p,q}(\phi_{\eps})} \pl .$$ This gives an
$(1+\eps)^{\frac2s}$-isomorphism
\begin{equation} \label{Ieps}
 \big[ L_p(N), \Delta_{p,q}(\phi) \big]_\theta = \big[ L_p(N),
 \Delta_{p,q}(\phi_{\eps}) \big]_\theta.
\end{equation}
In addition, $\Delta_{p,r}(\phi) = \Delta_{p,r}(\phi_\eps)$ are
$(1+\eps)^{\frac2u}$-isomorphic with $\frac1u = \frac1r - \frac1p$
and
 \[ \Delta_{p,r}(\phi) \lel \Delta_{p,r}(\phi_{\eps}) \lel \big[
 L_p(N), \Delta_{p,q}(\phi_{\eps}) \big]_\theta \lel \big[ L_p(N),
 \Delta_{p,q}(\phi) \big]_\theta \pl .\]
This proves the first assertion. Let us denote by $I_{\eps}:
\Delta_{p,r}(\phi_\eps) \to \Delta_{p,r}(\phi)$ the formal
identity. Let $\mathcal{Q}_r(\eps): L_r(N \oplus N) \to
\Delta_{p,r}(\phi_{\eps})$ be the projection constructed above.
Then we denote by $\mathcal{Q}_r: L_r(N \oplus N) \to
\Delta_{p,r}(\phi)$ the densely defined map
 \[ \mathcal{Q}_r (d^{\frac1u} \alpha, \beta d^{\frac1u}) \lel
 I_{\eps} \mathcal{Q}_r(\eps) (d_{\eps}^{\frac1u} \alpha,
 \beta d_{\eps}^{\frac1u}) \pl. \]
Since we have $\|\mathcal{Q}_r\| \le (1+\eps)^{\frac1r - \frac1p}
\|\mathcal{Q}_r(\eps)\|$, it turns out that $\mathcal{Q}_r$ is the
desired projection. \qd

\begin{rem}\label{sfform}{\rm Let us explain how we may pass
to the limit $\eps\to 0$ for the definition of $\mathcal{Q}_r$. We
denote by $B_{\infty}(\rz)$ the algebra of bounded measurable
functions on $\rz$ and find a normal $*$-representation
$\pi:B_{\infty}(\rz)\ten_{\min} B_{\infty}(\rz)\to B(L_2(N))$
given by $\pi(f\ten g)\lel \L_{f(d)}\R_{g(d)}$. This shows that
 \[ \L_{d^{1/u}}(\L_{d^{1/u}}+\R_{d^{1/u}})^{-1}
 \lel \mathrm{SOT}-\lim_{\eps\to 0}
 \L_{d_{\eps}^{1/u}}(\L_{d_{\eps}^{1/u}}+\R_{d_{\eps}^{1/u}})^{-1} \pl. \]
A similar statement holds for
$\R_{d^{1/u}}(\L_{d^{1/u}}+\R_{d^{1/u}})^{-1}$. Therefore, for
$x\in L_r(N)$ the family $T_{\eps}(x)=
\L_{d_{\eps}^{1/u}}(\L_{d_{\eps}^{1/u}}+\R_{d_{\eps}^{1/u}})^{-1}(x)$
is uniformly bounded in $L_r(N)$ and converges in $L_2(N)$. It
follows very easily from \cite[Theorem 3.6]{KF} that $T_{\eps}(x)$
converges in $L_r(N)$. We recall the canonical embedding $j_{p,r}:
\Delta_{p,r}(\phi) \to L_r(N\oplus N)$ given by $j_{p,r}(x) =
(\L_{d^{1/u}}x,\R_{d^{1/u}}x)$ and deduce that
 \[ j_{p,r}(\L_{d^{1/u}}+\R_{d^{1/u}})^{-1}: L_r(N \oplus N) \to
 L_r(N \oplus N) \]
is a well-defined bounded map. Thus $\mathcal{Q}_r \lel
(\L_{d^{1/u}}+\R_{d^{1/u}})^{-1}$ is a projection onto
$\Delta_{p,r}(\phi)$ and the pointwise limit of the
$\mathcal{Q}_r(\eps)$'s. In particular, the condition from Remark
\ref{propertyQ} is satisfied. Indeed, using the Borel functional
calculus for $B_{\infty}(\rz)\ten B_{\infty}(\rz)$ we find
 \[ \mathcal{Q}_r(y,z) \lel \int_{\rz \times \rz}
 (d(\om)^{1/u} + d(\om')^{1/u})^{-1}
 dE_{\om}(y+z)dE_{\om'} \pl .\]
Let us note that in the semifinite case (without assuming $c_1 \le
d \le c_2$ but still assuming $d$ is faithful), we may obtain the
same formula by using an increasing net of spectral projections of
$d$.}
\end{rem}

The proof for the general case is based on Haagerup's reduction
theorem, see \cite{jx-red}. Let us briefly explain how this
construction works. Let us consider a von Neumann algebra $N$
equipped with a normal faithful state $\phi$ associated to a
density $d$. Let us define the discrete group
$$\mathrm{G} = \bigcup_{n \in \mathbb{N}} 2^{-n} \mathbb{Z}.$$
Then we construct the crossed product $M = N \rtimes_{\sigma^\phi}
\mathrm{G}$. That is, if $H$ is the Hilbert space provided by the
GNS construction applied to $\phi$ and $\sigma^\phi$ denotes the
one parameter modular automorphism group on $N$ associated to
$\phi$, then $M$ is generated by the representations $\pi: N \to
\mathcal{B}(L_2(\mathrm{G}; H))$ and $\lambda: \mathrm{G} \to
\mathcal{B}(L_2(\mathrm{G}; H))$, where
$$\big( \pi(x) \xi \big) (g) = \sigma_{-g}^\phi(x) \xi(g) \qquad
\mbox{and} \qquad \big( \lambda(h) \xi \big) (g) = \xi(g-h).$$ By
the faithfulness of $\pi$ we are allowed to identify $N$ with its
image $\pi(N)$. Then, a generic element in the crossed product $M$
has the form $\sum_g x_g \lambda(g)$ with $x_g \in N$ and we have
the conditional expectation
 \[ \mathsf{E}_{N}\big(\sum_{g \in
 \mathrm{G}} x_g \lambda(g)\big) \lel  x_0 \in N \pl .\]
The algebra $M$ contains an increasing net $(M_\alpha)_{\alpha \in
\Lambda}$ of finite von Neumann subalgebras with normal
conditional expectations $\mathcal{E}_\alpha: M \to M_\alpha$. One
of the important properties of Haagerup's construction is that
$\psi = \phi \circ \mathsf{E}_N$ is a normal faithful state such
that $\psi \circ \mathcal{E}_\al = \psi$ holds for each $\alpha
\in \Lambda$. Moreover, the restriction $\psi_\alpha$ of $\psi$ to
$M_\alpha$ has a density $d_\alpha$ such that $$c_1(\alpha)
1_{M_\alpha} \le d_\alpha \le c_2(\alpha) 1_{M_\alpha}$$ for some
constants $0 < c_1(\alpha) \le c_2(\alpha) < \infty$. If $d_\psi$
denotes the density associated to the state $\psi$, we consider
the canonical conditional expectation $\mathcal{E}_{\alpha,p}:
L_p(M) \to L_p(M_\alpha)$ and the canonical inclusion
$\iota_{\alpha,p}: L_p(M_\alpha) \to L_p(M)$ densely defined
respectively by
$$\mathcal{E}_{\alpha,p}(x d_{\psi}^{\frac1p}) \lel
\mathcal{E}_\alpha(x) d_\alpha^{\frac1p} \quad \mbox{and} \quad
\iota_{\alpha,p}(x d_\alpha^{\frac1p}) \lel x d_{\psi}^{\frac1p}
\pl.$$ We refer to \cite{JX} for more information on these maps.
It is shown in \cite{jx-red} that
 \begin{equation} \label{k-limi}
 \limm_\alpha \iota_{\alpha,p} \hskip1pt \mathcal{E}_{\alpha,p}(x)
 \lel x \quad \mbox{for all} \quad x \in L_p(M) \quad \mbox{and}
 \quad 1 \le p < \infty \pl.
 \end{equation}
We will also need the $L_p$ version of $\mathsf{E}_N:L_p(M)\to
L_p(N)$:
 \[ \mathsf{E}_{N,p}(x d_{\psi}^{\frac1p})\lel
 \mathsf{E}_N(x) d^{\frac1p} \pl.\]
This comes with the natural inclusion map $j_{N,p}:L_p(N)\to
L_p(M)$, $j_{N,p}(x d^{\frac1p})=xd_{\psi}^{\frac1p}$, see again
\cite{JX}. With this information we start our approximation
procedure. Indeed, the following mappings will be instrumental in
our proof of Theorem \ref{minint} for general von Neumann algebras
$$u_{\alpha,p} = \mathsf{E}_{N,p} \hskip1pt \iota_{\alpha,p}:
L_p(M_\alpha) \to L_p(N) \quad \mbox{and} \quad w_{\alpha,p} =
\mathcal{E}_{\alpha,p} j_{N,p}: L_p(N) \to L_p(M_\alpha) \pl.$$

\begin{lemma}\label{refll}
The following properties hold:
\begin{itemize}
 \item[i)] If $1\le p<\infty$, $\lim_\alpha u_{\alpha,p}
 w_{\alpha,p} (x) = x$ for all $x\in L_p(N)$.

 \item[ii)] The mappings $u_{\alpha,p}$ and $w_{\alpha,p}$ induce
 contractions
 \begin{eqnarray*}
 u_{\alpha,p}: & \hskip1pt \big[ L_p(M_\alpha),
 \Delta_{p,q}(\psi_\alpha)\big]_{\theta} \to
 \big[ L_p(N), \Delta_{p,q}(\phi) \big]_{\theta} \pl, \\
  w_{\alpha,p}: & \big[ L_p(N), \Delta_{p,q}(\phi) \big]_{\theta}
 \to \big[ L_p(M_\alpha), \Delta_{p,q}(\psi_\alpha) \big]_{\theta} \pl.
 \end{eqnarray*}

 \item[iii)] If $1 \le q < p \le \infty$ and $0 < \theta < 1$, we have
 $$\limm_\alpha u_{\alpha,p} w_{\alpha,p}(x)=x \quad \mbox{for all}
 \quad x \in \big[ L_p(N), \Delta_{p,q}(\phi) \big]_{\theta}.$$
\end{itemize}
\end{lemma}

\begin{proof}
Since $\mathsf{E}_{N,p} j_{N,p} (x) = x$ for all $x \in L_p(N)$,
we have $$\limm_\alpha u_{\alpha,p} w_{\alpha,p} (x) - x =
\limm_\alpha \mathsf{E}_{N,p} \Big( \iota_{\alpha,p} \hskip1pt
\mathcal{E}_{\alpha,p} \big( j_{N,p}(x) \big) - j_{N,p}(x) \Big) =
0,$$ where the last identity follows from \eqref{k-limi} and the
contractivity of $\mathsf{E}_{N,p}$ in $L_p(M)$. This proves the
first assertion. Now let us identify $\Delta_{p,q}$ with its image
$j_{p,q}(\Delta_{p,q})$ in $L_q(N \oplus N)$ and also
$\Delta_{p,q}(\psi_\al)$ with its image
$j_{p,q}(\Delta_{p,q}(\psi_\al))$ in $L_q(M_\al \oplus M_\al)$.
Then, to prove ii) we will regard the mapping $$w_{\alpha,p}:
\Delta_{p,q}(\phi) \to \Delta_{p,q}(\psi_\alpha)$$ as the
restriction of $w_{\alpha,q} \oplus w_{\alpha,q}: L_q(N \oplus N)
\to L_q(M_\alpha \oplus M_\alpha)$ to the subspace
$$\Big\{ (d^{\frac1s} x, x d^{\frac1s}) \, \big| \ x \in L_p(N)
\Big\} \quad \mbox{with} \quad 1/s=1/q-1/p \pl .$$ If $x = y
d^{\frac1p}$ for $y \in N$, we have
\begin{align} \label{commwj}
 w_{\alpha,q} (j_{p,q}(x)) & \lel \big( w_{\alpha,q}
 (d^{\frac1s} x), w_{\alpha,q}(x d^{\frac1s}) \big) \\ \nonumber &
 \lel \big( d_\al^{\frac1s} \mathcal{E}_\alpha(y) \hskip1pt
 d_\alpha^{\frac1p}, \hskip2pt
 \mathcal{E}_\alpha(y) \hskip1pt
 d_\alpha^{\frac1q} \big) \lel  j_{p,q} (\mathcal{E}_{\alpha}(y)
d_\alpha^{\frac1p}) \lel j_{p,q}(w_{\alpha,p}(x)) \pl .
\end{align}
Here we use the well-known fact that
$\mathcal{E}_{\alpha,p}(d_{\psi}^{\frac{1-\eta}{p}} x
d_{\psi}^{\frac{\eta}{p}}) = d_\alpha^{\frac{1-\eta}{p}}
\mathcal{E}_\alpha(x) d_\alpha^{\frac{\eta}{p}}$, which follows
from our definition of $\mathcal{E}_{\alpha,p}$ and the identity
$\mathcal{E}_\alpha \sigma^{\psi} = \sigma^{\psi_\alpha}
\mathcal{E}_\alpha$,  see \cite{JX} for further details. Therefore
the map $w_{\alpha,p}$ induces a compatible contraction on the
interpolation couple $(L_p(N), \Delta_{p,q}(\phi))$ and hence on
the complex interpolation space $[L_p(N),
\Delta_{p,q}(\phi)]_{\theta}$. The argument for $u_{\alpha,p}$ is
entirely similar. In the proof of iii) we first observe that it
suffices to prove the assertion on a dense subspace, because we
already know from ii) that the maps $u_{\alpha,p} w_{\alpha,p}$
are contractions. If $x \in L_p(N)$ (we remind the reader that
$p=\infty$ is allowed and hence we may not assume that
$\lim_\alpha u_{\alpha,p} w_{\alpha,p}(x) = x$ holds in norm), we
set $\gamma_{\alpha,p} = u_{\alpha,p} w_{\alpha,p}$ and have
 \begin{eqnarray*}
 \lefteqn{\limm_\alpha \|\gamma_{\alpha,q} (x)
 - x\|_{[L_p(N), \Delta_{p,q}(\phi)]_{\theta}}}
 \\ & \le & \limm_\alpha \|\gamma_{\alpha,p} (x) -
 x\|_p^{1-\theta} \Big( \big\|
 d^{\frac1s} \big( \gamma_{\alpha,p}(x) - x \big) \big\|_q^q
 + \big\| \big( \gamma_{\alpha,p} (x) - x \big) d^{\frac1s}
 \big\|_q^q \Big)^{\frac{\theta}{q}} \\ & \le & (2 \hskip1pt
\|x\|_p)^{1-\theta} \limm_\alpha \Big( \big\| \gamma_{\alpha,q}
(d^{\frac1s} x) - d^{\frac1s} x \big\|_q^q + \big\|
\gamma_{\alpha,q} (x d^{\frac1s}) - x d^{\frac1s} \big\|_q^q
\Big)^{\frac{\theta}{q}} \lel 0 \pl.
\end{eqnarray*}
The first inequality uses the three lines lemma, the second
applies i) and uses $\theta>0$. \qd

\begin{proof}[{\bf Step 3 of the proof.}]

We now conclude the proof of Theorem \ref{minint}. For the
assertion a) we observe that the upper estimate in Remark
\ref{explicit} holds in general by the same argument used in Step
1 above. For the lower estimate we observe that $(M_\alpha,
d_\alpha)$ satisfies the hypotheses of Step 2. Hence we have
 \[ 
 \Big( \big\|
 d_\alpha^{\frac1r-\frac1p} x \big\|_r^r + \big\| x
 d_\alpha^{\frac1r-\frac1p} \big\|_r^r \Big)^{\frac1r} \sim
 \|x\|_{[L_p(M_\alpha),
 \Delta_{p,q}(\psi_\alpha)]_{\theta}} \pl \]
for all $x \in L_p(M_\alpha)$ and $\alpha \in \Lambda$. This
implies that
\begin{align*}
 \|x\|_{[L_p(N), \Delta_{p,q}(\phi)]_{\theta}} & \kl
 \limmsup_\alpha \|\gamma_{\alpha,p} \hskip3pt
 (x)\|_{[L_p(N), \Delta_{p,q}(\phi)]_{\theta}} \\
 & \kl \limmsup_\alpha \|w_{\alpha,p}(x)\|_{[L_p(M_\alpha),
 \Delta_{p,q}(\psi_\alpha)]_{\theta}} \\
 & \pl \lesssim  \ \,  \limmsup_\alpha
 \Big( \big\| d_\alpha^{\frac1r-\frac1p}
 \mathcal{E}_{\alpha,p}(x) \big\|_r^r + \big\|
 \mathcal{E}_{\alpha,p}(x) d_\alpha^{\frac1r-\frac1p}
 \big\|_r^r \Big)^{\frac1r} \\ & \kl \Big( \big\|
 d^{\frac1r-\frac1p} x \big\|_r^r + \big\| x
 d^{\frac1r-\frac1p} \big\|_r^r \Big)^{\frac1r} \pl .
\end{align*}
We will now construct the projection as a suitable limit. Let
 \[ \mathcal{Q}_{\alpha,r}: L_r(M_\alpha \oplus M_\alpha) \to
 \Delta_{p,r}(\psi_\alpha) \]
be the projection from Step 2 and let $\U$ be a free ultrafilter
on $\Lambda$. Then we define
 \[ \big\langle \mathcal{Q}_r(x,y), \xi \big\rangle
 \lel \limm_{\alpha,\U} \big\langle u_{\alpha,p}
 \mathcal{Q}_{\alpha,r} (w_{\alpha,r} (x), w_{\alpha,r}(y)), \xi
 \big\rangle \]
for every $\xi \in \Delta_{p,r}(\phi)^*$. Note that
$\Delta_{p,r}(\phi)$ is a reflexive Banach space. Therefore, we
deduce that we have $\mathcal{Q}_r(x,y) \in \Delta_{p,r}(\phi)$
for all $(x,y) \in L_r(N \oplus N)$. Since
$\mathcal{Q}_{\alpha,r}$ is a projection, we deduce
$$u_{\al,p}\mathcal{Q}_{\alpha,r} \big( w_{\alpha,r} (d^{\frac1r-\frac1p}
x), w_{\alpha,r}(x d^{\frac1r-\frac1p}) \big) \lel u_{\al,p}
\mathcal{Q}_{\alpha,r} \big( d_\alpha^{\frac1r-\frac1p}
w_{\alpha,p}(x), w_{\alpha,p}(x) d_\alpha^{\frac1r-\frac1p} \big)
= \gamma_{\alpha,p}(x) \pl .$$ Thus Lemma \ref{refll} iii) and
$[L_p(N), \Delta_{p,q}(\phi)]_{\theta} = \Delta_{p,r}(\phi)$ imply
 \[ \mathcal{Q}_r(d^{\frac1r-\frac1p} x, x d^{\frac1r-\frac1p}) =
 x \]
for all $x \in L_p(N)$. Since $\mathcal{Q}_r$ is continuous, we
deduce the result by density. \qd

\begin{rem}{\rm Theorem \ref{minint} also holds in the category
of operator spaces. That is, the map $\mathcal{Q}_r: L_r(N \oplus
N) \to \Delta_{p,r}(\phi)$ is completely bounded. This follows
immediately from replacing $d$ by $1\ten d$ in $L_1(M_m(N))$.
Moreover, in the semifinite setting the assumption $\tau(d)=1$ is
not really needed. Therefore, Theorem \ref{minint} also holds for
$\tau$-measurable operators $d$. More generally, this can be
extended to strictly semifinite weights.  At the time of this
writing it is not clear whether there is a result in this
direction for arbitrary weights. For two densities $d_1$ and $d_2$
we can obtain results in this direction by considering $(1,2)$
entries in the space $\Delta_{p,r}(\phi_2)$, where $\phi_2$ is
associated to the density $d = d_1 \ten e_{11} + d_2 \ten e_{22}$
on $M_2(N)$. We leave the details to the interested reader.}
\end{rem}

Using the methods of our paper, the referee found a proof for the
following interpolation result which generalizes our Theorem
\ref{minint}. We are indebted to the referee for allowing us to
reproduce his argument.

\begin{theorem}\label{pr}
Let $1 \le q_0, q_1 \le \infty$ and $\al_0,\al_1\gl 0$. Define
$1/q_{\theta} = (1-\theta)/q_0 + \theta/q_1$ and $\al_{\theta} =
(1-\theta)\al_0 + \theta \al_1$ for $0 < \theta < 1$. Then, the
following isomorphism holds for any density $d$ of a normal
faithful state on $N$ $$\Delta_{q_{\theta}} (N,d^{\al_{\theta}})
\lel \big[ \Delta_{q_0}(N,d^{\al_0}), \Delta_{q_1} (N,d^{\al_1})
\big]_{\theta} \pl .$$
\end{theorem}

\begin{proof} Here we will prove the result assuming
\eqref{simfun}. The proof in the general case follows by
approximation in the semifinite case and an application of
Haagerup's decomposition, as in Step 3 above. Using the triangular
map, it suffices to prove
 \[ \Delta_{q_{\theta}}^{UT}(N,d^{\al_{\theta}})
 = \big[ \Delta_{q_0}^{UT}(N,d^{\al_0}),
 \Delta_{q_1}^{UT}(N,d^{\al_1}) \big]_{\theta} \pl . \]
According to Lemma \ref{intpis0}, the direct sums
$UT_{q_{\theta}}^{r}(N,d^{\al_\theta}) \oplus
UT_{q_\theta}^{c}(N,d^{\al_\theta})$ are an interpolation scale.
Thus, it suffices to find a common projection which is bounded for
$q_0$ and $q_1$. Let us show that the map
 \[ \mathcal{Q}(y,z) \lel (x,x) \quad \mbox{where}\quad
  x\lel (\L_{d^{\al_0+\al_1}}+\R_{d^{\al_0+\al_1}})^{-1}
  \big(\L_{d^{\al_0+\al_1}}(y)+\R_{d^{\al_0+\al_1}}(z)\big) \]
is bounded in both spaces. Indeed, Lemma \ref{Schur} gives
\begin{eqnarray*}
\big\| (\L_{d^{\al_0+\al_1}}+\R_{d^{\al_0+\al_1}})^{-1}
\hskip1.5pt \L_{d^{\al_0+\al_1}}(y)
\big\|_{UT_{q_j}^r(N,d^{\al_j})}
& \le & \frac32 \pl \|y\|_{UT_{q_j}^r(N,d^{\al_j})} \pl , \\
\big\| (\L_{d^{\al_0+\al_1}}+\R_{d^{\al_0+\al_1}})^{-1}
\R_{d^{\al_0+\al_1}}(z) \big\|_{UT_{q_j}^c(N,d^{\al_j})} & \le &
\frac32 \pl \|z\|_{UT_{q_j}^c(N,d^{\al_j})} \pl ,
\end{eqnarray*}
for $j=0,1$. Hence, it remains to see that
\begin{eqnarray*}
\big\| (\L_{d^{\al_0+\al_1}}+\R_{d^{\al_0+\al_1}})^{-1}
\hskip1.5pt \L_{d^{\al_0+\al_1}}(y)
\big\|_{UT_{q_j}^c(N,d^{\al_j})}
& \le & \frac32 \pl \|y\|_{UT_{q_j}^r(N,d^{\al_j})} \pl , \\
\big\| (\L_{d^{\al_0+\al_1}}+\R_{d^{\al_0+\al_1}})^{-1}
\R_{d^{\al_0+\al_1}}(z) \big\|_{UT_{q_j}^r(N,d^{\al_j})} & \le &
\frac32 \pl \|z\|_{UT_{q_j}^c(N,d^{\al_j})} \pl .
\end{eqnarray*}
Since all these cross estimates can be handled similarly, we only
estimate the first one in the case $j =0$. Using $\eta = \alpha_0
/ (\alpha_0 + \alpha_1)$ in conjunction with Lemma 1.7 one more
time, we obtain
\begin{eqnarray*}
\lefteqn{\big\| (\L_{d^{\al_0+\al_1}}+\R_{d^{\al_0+\al_1}})^{-1}
\hskip1.5pt \L_{d^{\al_0+\al_1}}(y)
\big\|_{UT_{q_0}^c(N,d^{\al_0})}} \\ [4pt] & = & \big\|
\L_{d^{\al_1}} \R_{d^{\al_0}}
(\L_{d^{\al_0+\al_1}}+\R_{d^{\al_0+\al_1}})^{-1} \hskip1.5pt
\L_{d^{\al_0}}(y) \big\|_{UT_{q_0}(N)} \\ & \le & \frac32 \,
\|\L_{d^{\al_0}}(y)\|_{UT_{q_0}(N)} \lel \frac32 \, \|d^{\al_0}
y\|_{UT_{q_0}(N)} \lel \frac32 \, \|y\|_{UT_{q_0}^r(N,d^{\al_0})}
\pl .
\end{eqnarray*}
We apply the same arguments (and same ``projection'') for lower
triangular elements. \qd


\noindent The following byproduct of our arguments might be of
independent interest.

\begin{cor}\label{15}
Let $1 < p < \infty$ and $\alpha > 0$. Then the maps
\[ \R_{d^{(1-\eta)\al}}\L_{d^{\eta\al}}
 (\L_{d^{\al}}+\R_{d^{\al}})^{-1} \quad (0 \le \eta \le 1) \]
are bounded on $L_p(N)$ for any density $d$ of a normal faithful
state on $N$.
\end{cor}

\begin{proof} For $d=\sum_{k=1}^n d_k e_k$ as in \eqref{simfun}
this follows immediately from Lemma \ref{Schur} and Lemma
\ref{tri}. Then we follow the same procedure as in the proof of
Theorem \ref{minint} by first showing it for finite von Neumman
algebras with densities bounded above and below, and then apply
the Haagerup construction.\qd

\section{Subspaces of noncommutative $L_1$}

In this section we follow Pisier's approach and prove Theorem
\ref{main0}. Let us recall the notions of type and cotype from
Banach space theory. Given a probability space $\Omega$, let us
consider a sequence $(\eps_k)$ of independent Bernoulli random
variables equidistributed in $\pm 1$. A linear map $T:X\to Y$ has
type $p$ if there exists $c_1>0$ such that the inequality below
holds for all finite sequences $(x_k)$ in $X$
 \[ \Big( \mathbb{E} \,
 \big\| \summ_k \eps_k T(x_k) \big\|_Y^2 \Big)^{\frac12} \kl c_1 \Big(
  \summ_k \|x_k\|_X^p \Big)^{\frac1p} \pl. \]
Then $t_p(T)=\inf c_1$ satisfying the inequality above. A Banach
space has type $p$ if $id_X$ has type $p$. We use the standard
notation $t_p(X)=t_p(id_X)$. A linear map $T:X\to Y$ is said to be
of cotype $q$ if
 \[
 \Big( \summ_k \| T(x_k)\|_Y^q \Big)^{\frac1q} \kl c_2 \Big(
  \mathbb{E} \, \big\| \summ_k \eps_k x_k \big\|_X^2 \Big)^{\frac12}
  \pl. \]
We define $c_q(T)=\inf c_2$, where the infimum is taken over all
$c_2$ satisfying the inequality above. Again $c_q(X)=c_q(id_X)$
for a Banach space $X$.  Given a von Neumann algebra $N$, a linear
map $T: L_p(N)
 \to X$ is called $(q,+)$-summing if there exists a constant $c>0$
such that the inequality below holds for all finite sequence
$(x_k)$ of positive elements $x_k \in L_p(N)$
 \begin{equation} \label{q+summing}
 \Big( \summ_k \|T (x_k)\|_X^q \Big)^{\frac1q} \kl c \, \Big\|
 \summ_k x_k \Big\|_p \pl.
 \end{equation}
We denote $\pi_{q,+}(T)=\inf c$. Let us recall the well-known fact
 \begin{equation}
 \label{cotplus}
  \pi_{q,+}(T:L_p(N)\to X)\kl 2 c_q(T) \pl .
  \end{equation}
Indeed, for positive elements $x_k$ the order relation implies
that
 \[ \big\| \sum_k \eps_k x_k \big\|_p
 \kl \big\| \sum_{k,\eps_k=1}x_k \big\|_p + \big\|
 \sum_{k,\eps_k=-1}x_k \big\|_p
 \kl 2 \, \big\| \sum_k x_k \big\|_p \pl .\]
We shall also need the following well-known fact from
interpolation \cite[section 4.7]{BL}.

\begin{lemma} \label{fact} Let $(A_0,A_1)$ be an interpolation
couple of Banach spaces. Assume that $A_0$ is contractively
included in $A_1$ and let $0<\tilde{\eta}<\eta < 1$. Then, there
exists some absolute constant $c(\eta,\tilde{\eta})$ depending
only on $(\eta,\tilde{\eta})$ such that the norm of the inclusion
\[ [A_0,A_1]_{\tilde{\eta}} \subset [A_0,A_1]_{\eta,1} \quad \mbox{is
controlled by} \quad c(\eta,\tilde{\eta}) \pl.\]
\end{lemma}
\noindent Without assuming full support for $d$ we keep the
notation
\begin{equation} \label{notfullsupport}
\|x\|_{\Delta_{p,q}(\phi)} \lel \max \Big\{
\|d^{\frac1q-\frac1p}x\|_q, \|xd^{\frac1q-\frac1p}\|_q \Big\} \pl
.
\end{equation}
If we set $\mbox{supp} \, d = e$, the expression above vanishes on
$(1-e) L_p(N) (1-e)$. Relation \eqref{notfullsupport} defines a
norm on $e L_p(N) + L_p(N)e$, a complemented subspaces of the
quasi-normed space $(L_p(N),\|\pl\|_{\Delta_{p,q}})$. We will
write $\Delta_{p,q}(\phi)$ for the completion of $eL_p(N)+L_p(N)e$
with respect to this norm. The spaces $e \Delta_{p,q}(\phi) e$, $e
\Delta_{p,q}(\phi) (1-e)$ and $(1-e) \Delta_{p,q}(\phi) e$ are the
complemented subspaces of $\Delta_{p,q}(\phi)$ obtained from the
closure of $eL_p(N)e$, $eL_p(N)(1-e)$ and $(1-e)L_p(N)e$ in
$\Delta_{p,q}(\phi)$.

\begin{lemma}\label{three} Let $d$ be the density of a
normal state $\phi$ and let $e$ be the support projection of $d$,
so that $\phi$ is faithful on $eNe$. If
$\frac{1}{r}=\frac{1-\theta}{p}+\frac{\theta}{q}$, then $[eL_p(N)
+ L_p(N)e,\Delta_{p,q}(\phi)]_{\theta}$ is isomorphic to the
direct sum
 \[ e\Delta_{p,r}(\phi)e \oplus eL_r(N)(1-e) \oplus (1-e)L_r(N)e \pl .\]
The restriction of this isomorphism on $eL_p(N) + L_p(N)e$ is
given by
 \[ x \mapsto
 \Big( exe,d^{\frac1r-\frac1p}x(1-e), (1-e)xd^{\frac1r-\frac1p} \Big)\pl .\]
\end{lemma}

\begin{proof}
By definition, we have $$e \Delta_{p,q}(\phi) e =
\Delta_{p,q}(\phi_{|_{eNe}}) \pl .$$ Hence, we can apply Theorem
\ref{minint} and find that
 \[ [eL_p(N)e,e\Delta_{p,q}(\phi)e]_{\theta} \simeq
 e\Delta_{p,r}(\phi)e \pl .\]
Now we discuss the off-diagonal parts in
\begin{eqnarray}
\label{th}
 [eL_p(N)+L_p(N)e, \Delta_{p,q}(\phi)]_{\theta} & \simeq & [eL_p(N)e,
 e\Delta_{p,q}(\phi)e]_{\theta} \\ [+0.2cm] \nonumber & \oplus &
 [eL_p(N)(1-e),e\Delta_{p,q}(\phi)(1-e)]_{\theta} \\ [+0.2cm]
 \nonumber & \oplus &
 [(1-e)L_p(N)e,(1-e)\Delta_{p,q}(\phi)e]_{\theta} \pl.
\end{eqnarray}
However, for $x=ex(1-e)$ we have that
 \begin{equation} \label{th1}
  \|ex(1-e)\|_{\Delta_{p,q}(\phi)} = \max \Big\{ \|d^{\frac1q-\frac1p}ex(1-e)\|_q,
 \|ex(1-e)d^{\frac1q-\frac1p}\|_q \Big\} = \|d^{\frac1q-\frac1p}x(1-e)\|_q \pl .
 \end{equation}
A similar remark applies for $x=(1-e)xe$. Therefore, the
interpolation space simplifies considerably in the off-diagonal
terms. Applying Kosaki's interpolation theorem we formally obtain
 \begin{equation}
 \label{ke}
 \big[ d^{\frac1q-\frac1p}eL_p(N)(1-e), eL_q(N)(1-e)
 \big]_{\theta} \lel d^{\frac1q-\frac1r} L_r(N) (1-e) \pl .
 \end{equation}
However, $\phi$ does not have full support and we can not apply
Kosaki's theorem directly. Let $\psi_{1-e} = \lim_j \psi_j$ be a
strictly semifinite weight on $(1-e) N (1-e)$. Then $\psi =
\psi_{1-e} + \phi$ is a strictly semifinite weight on $N$. Let
$e_j \le 1-e$ be the support of $\psi_j$ (with associated density
$d_j$) and $f_j = e_j + e$. We may apply Kosaki's interpolation
theorem for $\phi_j = \psi_j + \phi$ and the sum of the commuting
densities $d_j + d$. Then we obtain
 \[ \big[ (d+d_j)^{\frac1q-\frac1p} L_p(f_j N f_j), L_q(f_j N f_j)
 \big]_{\theta} \lel (d+d_j)^{\frac1q-\frac1r} L_{r}(f_j N f_j)
 \pl.\]
Since the map $W(y) = e y (1-e)$ is a contraction on the spaces at
both sides above, we can replace $N$ by $eN(1-e)$ in the isometric
isomorphism since the resulting spaces are contractively
complemented. Thus we find
 \[ \big[ d^{\frac1q-\frac1p} L_p(N) e_j, eL_q(N)e_j
 \big]_{\theta} \lel d^{\frac1q-\frac1r} L_{r}(N) e_j \pl.\]
Passing to the limit for $j \to \infty$ yields \eqref{ke}. Note
that in  \eqref{th} and \eqref{ke} we used different topological
vector spaces for the interpolation couple $(A_0,A_1)$. In
\eqref{th}, $A_0+A_1=A_1$  is the completion of $e L_p(N) (1-e)$
with respect to the norm in $\Delta_{p,q}(\phi)$. On the other
hand, in \eqref{ke} we use $e L_q(N) (1-e)$ as the underlying
vector space. Since $d^{\frac1q-\frac1p}L_p(N)(1-e)$ is dense in
$eL_q(N)(1-e)$, we have calculated the interpolation space. Thus
\eqref{th}, \eqref{th1} and \eqref{ke} imply that
$$\|ex(1-e)\|_{[eL_p(N)(1-e),e\Delta_{p,q}(\phi)(1-e)]_{\theta}}
 \lel \|d^{\frac1r-\frac1p}x(1-e)\|_{L_r(N)} \pl .$$
Taking adjoints, we obtain the same conclusion for the space
$(1-e)L_p(N)e$. \qd

\begin{theorem} \label{p} Let $X$ be a Banach space and $T:X \to
L_1(N)$ be a linear map such that $T^*$ has cotype $q$. Then,
there exists a density $d$ on $N$ such that, for all $1 < p < q'$,
we may construct a bounded linear map $u: X \to L_p(N) $
satisfying
\[ T(x) \lel d^{1-\frac{1}{p}} u(x) + u(x) d^{1-\frac{1}{p}}
\quad \mbox{for all} \quad x\in X \pl . \] If moreover $X \subset
L_1(N)$ is a subspace, $u$ is an isomorphic embedding of $X$ into
$L_p(N)$.
\end{theorem}

\begin{proof} We recall from \eqref{cotplus} that $T^*:N \to X^*$ is
$(q,+)$-summing. Therefore, we deduce from Pisier's factorization
theorem \cite[Theorem 3.2]{P1} that there exists a state $\phi$ on
$N$ such that
 \[ \|T^*(y)\|_{X^*} \kl c \, \|y\|^{1-\frac{2}{q}} \, \big(
 \phi(yy^*) + \phi(y^*y) \big)^{\frac{1}{q}}  \pl .\]
We use a standard trick (see \cite{P4}) to replace $\phi$ by its
normal part. Let $\phi_n$ be the normal part of $\phi$. Let
$(s_{\al})$ be a net of contractions in $N$ such that $\lim_{\al}
s_\al = 1$ in the strong operator and the strong$^*$ topology and
$\lim_{\al} \phi(s_{\al}ys_{\al})=\phi_n(y)$. Let $y\in N$ and
$x\in X$ of norm $1$ such that
 \[ \|T^*(y)\|_{X^*} \le (1+\eps) |T^*(y)(x)|\lel (1+\eps)  |tr(yT(x))| \pl .\]
We may write $T(x)=a=a_1a_2$ with $a_i\in L_2(N)$. Note that
 \[
 \limm_{\al} tr(s_{\al}ys_{\al}a) \lel \limm_\alpha
 tr(ys_{\al}as_{\al}) \lel tr(ya) \]
because
$\lim_{\al}s_{\al}as_{\al}=\lim_{\al}s_{\al}a_1a_2s_{\al}=a_1a_2=a$.
Therefore we find that
\begin{eqnarray*}
 \|T^*(y)\|_{X^*} & \le & (1+\eps) \limm_{\al}
 |tr(s_{\al}ys_{\al}a)|
 \\ & \le & c \, \limmsup_{\al} \|s_{\al}ys_{\al}\|^{1-\frac2q} \,
 \limm_{\al} \big( \phi(s_{\al}ys_{\al}y^*s_{\al}) +
 \phi(s_{\al}y^*s_{\al} y s_{\al}) \big)^{\frac{1}{q}} \\ & \le & c
 \, \limmsup_{\al} \|s_{\al}ys_{\al}\|^{1-\frac2q} \, \limm_{\al}
 \big( \phi(s_{\al}yy^*s_{\al}) + \phi(s_{\al}y^*y s_{\al})
 \big)^{\frac{1}{q}} \\ & \le & c \, \|y\|^{1-\frac2q} \big(
 \phi_n(yy^*) + \phi_n(y^*y) \big)^{\frac{1}{q}} \pl .
\end{eqnarray*}
Therefore, we may assume with no loss of generality that the state
$\phi$ is normal. This means that it is given by $\phi(y) =
tr(dy)$ for some density $d \in L_1(N)$. Let $e$ be the support of
$d$, so that $\phi$ is faithful on $eNe$. We then have
\begin{equation} \label{control}
 \|T^*(y)\|_{X^*} \kl c \, \|y\|_{N}^{1-\frac2q} \, \max \Big\{
 \|d^{\frac12} y\|_2, \, \|y d^{\frac12}\|_2 \Big\}^{2/q} \lel c \,
 \|y\|_N^{1-\frac2q} \, \|y\|_{\Delta_{2}(\phi)}^{2/q} \pl.
\end{equation}
Note that $T^*((1-e)y(1-e))=0$. According to a well-known result
(\emph{cf}. \cite[p.49]{BL}) we have
\begin{equation} \label{normal}
 \|T^*(y)\|_{X^*} \kl c \, \|y\|_{[eN+Ne,
 \Delta_{2}(\phi)]_{\frac2q,1}} \quad \mbox{for all} \quad y
 \in N \pl.
\end{equation}
Recall that we write $\Delta_q(\phi)$ for
$\Delta_{\infty,q}(\phi)$. We consider $(\eta,\tilde{\eta}) =
(2/q,2/p')$ and observe that $0 < \tilde{\eta} < \eta < 1$ since
$1 < p < q'$ and $2 < q < \infty$ (if $T^*$ has cotype $2$ it also
has cotype $q$ for all $q>2$). According to Lemma \ref{fact} we
deduce that
 \begin{equation} \label{normal2}
 \|T^*(y)\|_{X^*} \kl c \, c(\eta,\tilde{\eta})\,
 \|y\|_{[eN+Ne, \Delta_{2}(\phi)]_{\frac{2}{p'}}} \quad \mbox{for all} \quad y
 \in N \pl.
\end{equation}
Therefore, the map $T^*$ extends to a bounded map $T^*_{p'}:
[eN+Ne, \Delta_{2}(\phi)]_{2/p'}\to X^*$. Let us note that, in
accordance with Lemma \ref{three}, the intersection in this
interpolation space is $eNe+eN(1-e)+(1-e)Ne$ while by
\eqref{control} the map $T^*$ vanishes on the remaining corner
$(1-e)N(1-e)$. Let us recall the projection given by Theorem
\ref{minint}
\[ \mathcal{Q}_{p'}: L_{p'}(eNe \oplus eNe) \to
  e\Delta_{p'}(\phi)e \simeq \big[ eNe,
  e\Delta_{2}(\phi)e
  \big]_{\tilde{\eta}} \pl .\]
Using Lemma \ref{three} we may define the  map
$\widetilde{\mathcal{Q}}_{p'}: L_{p'}(N \oplus N) \to
\Delta_{p'}(\phi)$ by
 \[ \widetilde{\mathcal{Q}}_{p'}(y_1,y_2) \lel
 \mathcal{Q}_{p'}(ey_1e,ey_2e) \oplus ey_1(1-e) \oplus
 (1-e)y_2e\pl .\]
Thus by construction we have
 \begin{equation}\label{stt1}
  T^*(y) \lel  T_{p'}^*
  \widetilde{\mathcal{Q}}_{p'} (d^{\frac{1}{p'}} y,
  y d^{\frac{1}{p'}}) \end{equation}
for all $y \in N$. Unfortunately, $\widetilde{\mathcal{Q}}_{p'}$
does not vanish on vectors of the form $(y,-y)$. Therefore we need
a slight modification of $(\widetilde{\mathcal{Q}}_{p'})^*$ which
allows us to construct one map $u$ as asserted. For this we define
a map $v: \Delta_{p'}(\phi)^*\to L_p(N)$ as follows. According to
Lemma \ref{three} we have
$\Delta_{p'}(\phi)^*=(e\Delta_{p'}(\phi)e)^*\oplus
(eL_{p'}(N)(1-e))^*\oplus (1-eL_{p'}(N)e^*$. Following Remark
\ref{propertyQ} we know that
$\mathcal{Q}_{p'}^*(\xi)=(v_e(\xi),v_e(\xi))$ holds for some
bounded linear map $$v_e:(e\Delta_{p'}(\phi)e)^*\to L_{p}(eNe).$$
Thus we may define
 \[ v(\xi) \lel v_e(\xi_{e})+ \xi_{e,1-e}+\xi_{1-e,e} \quad
 \mbox{where $\xi$ has the components} \quad
 \xi=(\xi_e,\xi_{e,1-e},\xi_{1-e,e}) \pl .\] Under the usual
 duality bracket $\langle
a,b\rangle=tr(ab)$, we observe that
 \[ \big( eL_{p'}(N)(1-e) \big)^*=(1-e)L_p(N)e \quad \mbox{and} \quad
  \big( (1-e)L_{p'}(N)e \big)^*=eL_p(N)(1-e) \pl . \]
Therefore, we may and will assume that
$\xi_{e,1-e}=(1-e)\xi_{e,1-e}e$ and
$\xi_{1-e,e}=e\xi_{1-e,e}(1-e)$ are elements in $L_p(N)$. Then we
observe that
 \begin{align*}
  d^{\frac{1}{p'}}v(\xi)+v(\xi)d^{\frac{1}{p'}}
  &=d^{\frac{1}{p'}}v_{e}(\xi_e)+
   v_{e}(\xi_e)d^{\frac{1}{p'}} +
  d^{\frac{1}{p'}}\xi_{1-e,e}+\xi_{e,1-e}d^{\frac{1}{p'}} \pl .
  \end{align*}
This implies that, for all $y \in N$, we have
\begin{eqnarray*}
\lefteqn{\widetilde{\mathcal{Q}}_{p'}(d^{\frac{1}{p'}}y,yd^{\frac{1}{p'}}),\xi
\rangle} \\ & = & \langle
\widetilde{Q}_{p'}(d^{\frac{1}{p'}}eye,eye
d^{\frac{1}{p'}}),\xi_e\rangle +
 \langle d^{\frac{1}{p'}}y(1-e),\xi_{e,1-e}\rangle +
 \langle (1-e)yd^{\frac{1}{p'}},\xi_{1-e,e}\rangle \\
 & = &  tr(v_e(\xi_e)(d^{\frac{1}{p'}}eye+eyed^{\frac{1}{p'}}))  +
 tr(\xi_{e,1-e}d^{\frac{1}{p'}}y(1-e)) +
 tr(\xi_{1-e,e}(1-e)yd^{\frac{1}{p'}}) \\
 & = & tr((d^{\frac{1}{p'}}v(\xi)+v(\xi)d^{\frac{1}{p'}})y).
 \end{eqnarray*}
This will allow us to conclude easily. Indeed, we define
$u=v(T_{p'}^*)^*: X\to L_p(N)$. Then we deduce from \eqref{stt1}
that
 \begin{eqnarray*}
  tr((d^{\frac{1}{p'}}u(x)+u(x)d^{\frac{1}{p'}})y)
  & = & \big\langle
  \widetilde{\mathcal{Q}}_{p'}(d^{\frac{1}{p'}}y,
  yd^{\frac{1}{p'}}),(T_{p'}^*)^*(x) \big\rangle \\
  & = & \big\langle
  T_{p'}^*\widetilde{\mathcal{Q}}_{p'}(d^{\frac{1}{p'}}y,
  yd^{\frac{1}{p'}}),x \big\rangle \lel
  \big\langle T^*(y),x \big\rangle \lel tr(yT(x))
  \end{eqnarray*}
holds for all $y$. This means
$T(x)=d^{\frac{1}{p'}}u(x)+u(x)d^{\frac{1}{p'}}$. Let us now
consider the special case $T=\iota_X:X\to L_1(N)$ such that $X^*$
has cotype $q$. Then the left inverse for $u$ is given by $v(x) =
d^{\frac{1}{p'}} x + x d^{\frac{1}{p'}}$. Clearly, $v$ is bounded
and $u$ becomes an isomorphism. \qd

\begin{rem}{\rm The proof shows that we can construct the same $u$
because we only care about the restriction of $\mathcal{Q}_{p'}$
to elements of the form $(d^{1/p'}y,yd^{1/p'})$. If $N$ is
$\si$-finite and semifinite, we may assume $e=1$ and Remark
\ref{sfform} gives
 \[ u(x)\lel \int_{\rz\times \rz}
 [d(s)^{1/p'}+d(t)^{1/p'}]^{-1} dE_sT(x)dE_t \pl . \]
}\end{rem}
\begin{proof}[{\bf Proof of Theorem \ref{main0}.}]
The type index of $X$ is defined as
$$p_X = \inf \big\{ p \, | \, X \ \mbox{has type} \ p \big\}.$$
According to the Krivine-Maurey-Pisier theorem (see e.g. \cite{MS}
or \cite{Psv}) we know that for $p = p_X$ the spaces $\ell_p(n)$
are uniformly contained in $X$. If $p_X = 1$, we know from
\cite[Theorem 5.1]{RayXu} that $X$ contains $\ell_1$. However,
this contradicts the reflexivity of $X$. Hence, $p_X$ must be
strictly bigger that $1$. Let $p_0>1$ such that $X$ has type
$p_0$. This implies that $X^*$ has (finite) cotype $p_0'$ and
therefore Theorem \ref{p} applies. \qd

\begin{proof}[{\bf Proof of Corollary \ref{subseq}}]
Let $(x_n)$ be subsymmetric in $N_*$ and let $$X = {\rm span}
\big\{ x_n \, | \, n \ge 1 \big\} \pl .$$ According to (the proof
of) Theorem \ref{main0}, if $X$ does not contain $\ell_1$ then $X$
is isomorphic to a subspace of $L_p(N)$ for some $1<p<2$. Since we
know from \cite{JR} that $L_p(N)$ is asymptotically symmetric, we
deduce that $(x_n)$ is indeed symmetric. \qd

\begin{rem}{\rm Let $(x_n)$ be a subsymmetric
sequence in $L_1(N)$. A close inspection of \cite[Proposition
5.3]{RayXu} shows that $(x_n)$ contains a subsequence equivalent
to the unit vector basis of $\ell_1$ or $(x_n)$ is
$1$-equiintegrable (equivalently relatively weakly compact).
However, a subsymmetric sequence is equivalent to every
subsequence. Thus either $(x_n)$ is equivalent to the $\ell_1$
basis (hence symmetric) or $1$-equiintegrable. Therefore, the only
possibility of a subsymmetric, not symmetric sequence, occurs for
$1$-equiintegrable sequences where the unit ball of ${\rm
span}\{x_n:\nen\}$ is not $1$-equiintegrable, see also
\cite[Theorem 5.1]{RayXu}. }\end{rem}

\section{Nikishin-type results for $p$ finite}

In the commutative setting, Nikishin type results can be obtained
from a careful analysis of the maximal function. Although maximal
functions have been recently introduced in the noncommutative
setting \cite{JD,JX3}, they seem not to be applicable for this
type of results. Our approach using duality in the noncommutative
setting reduces the problem to norm estimates for positive
operators. In this section we prove the differential inequality
\eqref{diff00} and Theorem \ref{mainq}. Let us start with an
elementary observation. The result is known due to the work of
Araki \cite{Ar1} and Kosaki \cite{Kos3}. We give a short proof to
keep the paper more self-contained.

\begin{lemma} \label{ab}
Let $1 \le q \le \infty$ and $\al,\beta$ be positive. Then
 \[ \|\al^{\eta} \beta^{\eta}\|_{\frac{q}{\eta}} \le \|\al
 \beta\|_q^{\eta} \quad \mbox{for all} \quad 0 < \eta < 1 \pl .\]
\end{lemma}

\begin{proof}
Let us first show this for $\eta=\frac12$. Indeed,
 \[ \|\al^{\frac12} \beta^{\frac12}\|_{2q}^2 \lel
 \|\al^{\frac12} \beta \al^{\frac12}\|_q \pl .\]
Define $f(z) = \al^{1-z} \beta \al^{z}$ and fix $\lambda = \|\al
\beta\|_q$. We clearly have
 \[ \max \Big\{ \sup_{z \in \partial_0} \|f(z)\|_q, \, \sup_{z \in
 \partial_1} \|f(z)\|_q \Big\} \le \lambda \pl .\]
Therefore, we have $\|f(\frac12)\|_q \le \lambda$ and deduce the
assertion for $\eta = 1/2$. Now we show the inequality for all
$\frac12 < \eta < 1$. Take $c \in L_{(\frac{q}{\eta})'}(N)$ of
norm less than $1$. We may write $\eta/q = (1-\theta)/q +
\theta/2q$ for some $0 < \theta < 1$. Now we use interpolation and
assume that $N$ is $\sigma$-finite. The general case follows from
a well-known approximation argument. Using Kosaki's interpolation
theorem, we find an analytic function $g: \mathcal{S} \to
L_{(2q)'}(N)$ such that $g(\theta)=c$ and
 \[ \max \Big\{ \sup_{z \in \partial_0} \|g(z)\|_{q'}, \, \sup_{z
 \in \partial_1} \|g(z)\|_{(2q)'} \Big\} \le 1 \pl .\]
Therefore, the function
 \[ h(z) \lel tr \big( g(z) \al^{1-\frac{z}{2}}
 \beta^{1-\frac{z}{2}} \big) \]
is analytic. Here $tr$ denotes the trace on the Haageup $L_1$ space.
By the three lines lemma, we find
 \[ |tr(c \al^{\eta} \beta^{\eta})| \lel  |h(\theta)| \kl
 \big( \sup_{z \in \partial_0} |h(z)| \big)^{1-\theta} \big(
 \sup_{z \in \partial_1} |h(z)| \big)^\theta \pl .\]
However, we have
 \[ \sup_{z \in \partial_0} |h(z)| \kl \sup_{z \in
 \partial_0} \|g(z)\|_{q'} \|\al^{-z/2} \al \beta \beta^{-z/2}\|_q
  \kl \lambda \pl,\]  and
 \[ \sup_{z \in \partial_1} |h(z)| \kl \sup_{z
 \in \partial_1} \|g(z)\|_{(2q)'} \|\al^{-\mathrm{Im}(z)/2}
 \al^{\frac12} \beta^{\frac12} \beta^{-\mathrm{Im}(z)/2}\|_{2q} \kl
 \sqrt{\lambda} \pl.\]
Hence $|tr(c \al^{\eta} \beta^{\eta})| \le \lambda^{1-\theta}
\lambda^{\theta/2} = \lambda^{\eta}$. Finally, we observe that our
first argument for $\eta=1/2$ shows that if $\eta$ satisfies the
assertion, then so does $\eta/2$. Since the assertion holds for
$1/2\le \eta\le 1$, this completes the proof. \qd

\begin{theorem} \label{diff}
If $2< p<\infty$ and $a,x \in L_p(N)_+$, we have
 \[ \|a+x\|_p^p-\|a\|_p^p \kl p \, 2^{p-1} \max \Big\{
 \|a^{p-1}x\|_1, \|x\|_p^p \Big\} \pl .\]
\end{theorem}

\begin{proof}
We begin by recalling Lemma 3.1 (part 1) of \cite{Kos2}. In this
paper Kosaki used the uniform smoothness of $L_p(N)$ to show that
the function $f(s) = tr((a+sx)^p)$ is differentiable with
derivative
 \[ f'(s)\lel p \, tr \big( (a+sx)^{p-1}x \big) \pl .\]
This gives
 \begin{equation}\label{diffeqt}
  tr((a+x)^p)-tr(a^p) \lel p \int_0^1 tr \big( (a+sx)^{p-1}x \big)
  \, ds \pl .
  \end{equation}
We define $k$ to be the natural number satisfying $k\le p-1<k+1$ and
define
 \[ \theta \lel  \frac{p-1-k}{p-1} \quad \mbox{and} \quad
 \Big( \frac1q, \frac1r \Big) = \Big( \frac{k + 1-\theta}{p} \, , \,
 \frac{p-k-1+\theta}{p} \Big)\pl .\]
This implies $1/r=\theta$ and $1/q=1-\theta$. Then we may use
H\"older's inequality and find
\begin{eqnarray} \label{Klder}
tr \big( (a+sx)^{p-1} x \big) & = & tr \big( (a+sx)^{k}
x^{1-\theta} x^{\theta} (a+sx)^{p-1-k} \big) \\ \nonumber & \le &
\big\| (a+sx)^{k} x^{1-\theta} \big\|_q \pl \big\| x^{\theta}
(a+sx)^{p-1-k} \big\|_r \pl .
\end{eqnarray}
By Lemma \ref{ab} for $(\al,\beta,\eta) = (x^{\frac{1}{p-1}},
a+sx, p-1-k)$, we get
 \[ \big\| x^{\theta} (a+sx)^{p-1-k} \big\|_r \kl
 \big\| x^{\frac{1}{p-1}} (a+sx) \big\|_{p-1}^{p-1-k} \pl .\]
On the other hand, Lemma \ref{ab} for $(\al,\beta,\eta) = (x,
a^{p-1}, 1/(p-1))$ gives
\begin{align*}
 \big\| x^{\frac{1}{p-1}} (a+sx) \big\|_{p-1} & \kl
 \|x^{\frac{1}{p-1}}a\|_{p-1} + \|x^{\frac{p}{p-1}}\|_{p-1}
 \\ & \kl \|xa^{p-1}\|_1^{\frac{1}{p-1}} + \|x\|_p^{\frac{p}{p-1}}
 \kl  2 \, \max \Big\{ \|a^{p-1}x\|_1,\|x\|_p^p
 \Big\}^{\frac{1}{p-1}} \pl.
\end{align*}
Let us set
 \[  \xi \lel  \max \Big\{ \|a^{p-1}x\|_1,\|x\|_p^p \Big\}\pl .\]
Then we find the following estimate for the second term on the
right of \eqref{Klder}
\begin{equation} \label{00a}
\big\| x^{\theta} (a+sx)^{p-1-k} \big\|_r \kl  2^{p-1-k}
\xi^{\theta} \pl.
\end{equation}
We now consider the first term. For a subset $A \subset \{1,2,
\ldots, k\}$ we use the notation $a_{A^c}x_A=y_1\cdots y_k$ where
$y_i=x$ if $i \in A$ and $y_i=a$ if $i \in A^c$. Then we deduce
from the triangle inequality that
$$\big\| (a+sx)^{k} x^{1-\theta} \big\|_q \kl \summ_{A}
s^{|A|} \|a_{A^c} x_A x^{1-\theta}\|_q \kl \summ_{A} \|a_{A^c} x_A
x^{1-\theta}\|_q.$$ We claim that
 \begin{equation} \label{holdd}
 \|a_{A^c}x_Ax^{1-\theta}\|_{q} \kl
 \|a^{p-1}x\|_1^{\frac{k-|A|}{p-1}} \pl
 \|x\|_p^{1-\theta+|A|-\frac{k-|A|}{p-1}} \pl .
 \end{equation}
Before proving our claim, let us show how to finish the argument
 \begin{align*}
 & \summ_A \|a_{A^c} x_A x^{1-\theta}\|_q  \lel \sum_{j=0}^k
  {k\choose j} \, \|a^{p-1}x\|_1^{\frac{k-j}{p-1}} \|x\|_p^{1-\theta
 + j-\frac{k-j}{p-1}} \\
  & \le   \sum_{j=0}^k {k\choose j} \,
 \xi^{\frac{(k-j)}{p-1} + \frac{1-\theta}{p} + \frac{j}{p}
 -\frac{k-j}{p(p-1)}} \lel  \sum_{j=0}^k {k\choose j} \,
 \xi^{\frac{k}{p-1} + \frac{1-\theta}{p} -\frac{k}{p(p-1)}} \lel
  2^{k} \, \xi^{\frac{k+1-\theta}{p}} \lel 2^k \xi^ {1-\theta} \pl,
 \end{align*}
where the last identity follows from $1-\theta=\frac{k}{p-1}$. The
assertion then follows from the combination of (\ref{diffeqt},
\ref{Klder}, \ref{00a}) with the estimate given above. Therefore,
it remains to prove our claim. We need to consider different
cases. First assume that $A = \emptyset$, so that we have to prove
(recall that $1-\theta = \frac{k}{p-1}$) the inequality $\|a^k
x^{1-\theta}\|_{q} \le \|a^{p-1}x\|_1^{1-\theta}$. This follows
from Lemma \ref{ab} applied to $(\al,\beta,\eta) = (a^{p-1}, x,
1-\theta)$. Now assume $|A|\ge 1$. Then we may write
 \[ a_{A^c}x_A\lel a^{\al_1}x^{\beta_1}a^{\al_2}\cdots
 x^{\beta_r}a^{\al_{r+1}} \]
where $\sum_i \al_i+\sum_i \beta_i=k\le p-1$. Since we have
excluded the case $A = \emptyset$, all the coefficients $\al_i,
\beta_i$ are strictly positive, except possibly $\al_1$ and
$\al_{r+1}$. Let us first consider the case $\al_1>0=\al_{r+1}$.
We define $q_j$ for $1 \le j \le r$ by $1/q_j = (1+\al_j)/p$. Note
that $1 \le q_j \le p$ for all $j$. Then we use H\"{o}lder's
inequality and get
 \[ \|a_{A^c}x_Ax^{1-\theta}\|_{q} \kl
 \Big( \prod_{j=1}^r \|a^{\al_j}x\|_{q_j} \|x\|_p^{\beta_j-1} \Big)
 \|x\|_p^{1-\theta} \pl .\]
By Kosaki's interpolation theorem, we may estimate
 \[ \|a^{\al_j}x\|_{q_j} \kl \|x\|_p^{1-\theta_j}
 \|a^{p-1}x\|_1^{\theta_j} \]
where $\frac{1}{q_j}=\frac{1-\theta_j}{p}+\frac{\theta_j}{1}$.
This means $\theta_j=\frac{\al_j}{p-1}$. Therefore we find
$$\|a_{A^c}x_Ax^{1-\theta}\|_{q} \kl \|x\|_p^{1-\theta}
\prod_{j=1}^r \|a^{p-1}x\|_1^{\frac{\al_j}{p-1}}
\|x\|_p^{\beta_j-\frac{\al_j}{p-1}} \lel
\|a^{p-1}x\|_1^{\frac{k-|A|}{p-1}} \|x\|_p^{1-\theta +
|A|-\frac{k-|A|}{p-1}} \pl.$$ This proves \eqref{holdd} for
$\alpha_1 > 0 = \alpha_{r+1}$. Let us now also assume that
$\al_1=0$. Then we define the index $\widetilde{q}$ by
$1/\widetilde{q} = (\beta_1+ \al_2 + \beta_2)/p$. This allows us
to apply H\"older's inequality as above and obtain
 \[  \|a_{A^c}x_Ax^{1-\theta}\|_{q}
 \kl \|x^{\beta_1} a^{\al_2} x^{\beta_2} \|_{\widetilde{q}} \,
 \Big( \prod_{j=3}^r \|a^{\al_j} x \|_{q_j}
 \|x\|_p^{\beta_j-1} \Big) \|x\|_p^{1-\theta} \pl .\]
We can assume without loss of generality (taking adjoints if
necessary) that $\beta_1 \le \beta_2$. Define the index
$\widehat{q}$ by $1/\widehat{q} = (2\beta_1 + \al_2)/p$. Then we
deduce the following estimate from H\"older's inequality and Lemma
\ref{ab} applied to $(\al,\beta,\eta) = (x^{2\beta_1}, a^{\al_2},
\frac12)$
 \begin{align*}
 \|x^{\beta_1} a^{\al_2} x^{\beta_2}\|_{\widetilde{q}} & \kl
  \|x^{\beta_1} a^{\al_2} x^{\beta_1}\|_{\widehat{q}} \,
  \|x\|_p^{\beta_2-\beta_1}  \lel  \|x^{\beta_1}
  a^{\frac{\al_2}{2}}\|_{2 \widehat{q}}^2 \,
  \|x\|_p^{\beta_2-\beta_1} \\
 & \kl   \|x^{2\beta_1} a^{\al_2}\|_{\widehat{q}} \
 \|x\|_p^{\beta_2-\beta_1} \kl  \|xa^{\al_2}\|_{q_2} \,
  \|x\|_p^{\beta_2+\beta_1-1} \lel \prod_{1 \le j \le 2}
 \|a^{\alpha_j} x\|_{q_j} \|x\|_p^{\beta_j-1} \pl .
 \end{align*}
Therefore, the argument from above yields \eqref{holdd} in this
case. Thus we have treated the cases $\al_1=0=\al_{r+1}$ and
$\al_1>0=\al_{r+1}$. If $\al_1=0<\al_{r+1}$, we can take adjoints
and use the same argument one more time. Let us now assume
$\al_1>0$ and $\al_{r+1}>0$. If $\beta_r\ge 2$ the argument above
applies by splitting
$a^{\al_r}x^{\beta_r}a^{\al_{r+1}}=(a^{\al_r}x)x^{\beta_r-2}
(xa^{\al_{r+1}})$. Thus, the only case not covered so far is
$\beta_r=1$. Here we have to use a little trick
 \[  \|a^{\al_r}xa^{\al_{r+1}}\|_{\check{q}} \kl \|a^{\al_r +
 \al_{r+1}}x\|_{\check{q}} \quad \mbox{for} \quad 1/\check{q} =
 (\al_r + \al_{r+1} + 1)/p \pl .\]
Indeed, we define $d_r = a^{\al_r+\al_{r+1}}$ and $\gamma =
\frac{\al_r}{\al_r+\al_{r+1}}$. Then
 \[ \|a^{\al_r}xa^{\al_{r+1}}\|_{\check{q}} \lel \|d_r^{\gamma} x
 d_r^{1-\gamma}\|_{\check{q}} \pl .\]
Since the index $\check{q} \ge 1$, we may use complex
interpolation and define the analytic function $f(z) = d_r^z x
d_r^{1-z}$ on the strip. Then, the three lines lemma combined with
the fact that $x$ is self-adjoint implies that
 \[ \|f(\gamma)\|_{\check{q}} \le \max \Big\{ \sup_{z \in \partial_0}
 \|d_r^{z} x d_r^{1-z}\|_{\check{q}} \, , \, \sup_{z \in \partial_1}
 \|d_r^{z} x d_r^{1-z}\|_{\check{q}} \Big\}
 \le \max \Big\{ \|xd_r\|_{\check{q}}, \|d_rx\|_{\check{q}} \Big\} =
 \|d_r x\|_{\check{q}} \pl .\]
This allows us to repeat the same argument and thereby completes
the proof of \eqref{holdd}. \qd

\begin{rem} {\rm If $N$ is commutative,
the triangle inequality gives
 \begin{eqnarray*}
 \|a+x\|_p^p-\|a\|_p^p & = & p \int_0^1 tr \big( (a+sx)^{p-1}x
 \big) \pl ds \\ & \le & p \int_0^1 \Big(
 tr(a^{p-1}x)^{\frac{1}{p-1}} + s tr(x^p)^{\frac{1}{p-1}}
 \Big)^{p-1}ds \\ [4pt] & \le & (2^p-1) \, \max \Big\{ tr(a^{p-1}x),
 tr(x^p) \Big\} \pl.
 \end{eqnarray*}
However, in the noncommutative case it is known that the
expression $\phi(|x|^q)^{1/q}$ does not define a norm for
arbitrary states. On the other hand, for $2\le p\le 3$ Theorem
\ref{diff} follows immediately from the fact that $t\mapsto
t^{p-1}$ is operator convex. Indeed, we have
 \[ (a+sx)^{p-1} \lel (1+s)^{p-1} \Big( \frac{1}{1+s}a +
 \frac{s}{1+s}x \Big)^{p-1}
 \kl (1+s)^{p-2} (a^{p-1}+sx^{p-1}) \pl .\]
This implies
 \begin{align*}
 p \int_0^1 tr \big( (a+sx)^{p-1}x \big) ds
 &\le \frac{p(2^{p-1}-1)}{p-1} \big( tr(a^{p-1}x)+tr(x^p) \big)
 \pl .
 \end{align*}
This is even better than our estimate. }
\end{rem}

\begin{lemma} \label{pos}
Let $d$ be the density of a normal state and consider the norm
 \[ \|x\|_{p,t,d} \lel \max \Big\{ t^{\frac1p} \|x\|_p, \, t
 \|d^{\frac{1}{p'}}x\|_1, \, t \|xd^{\frac{1}{p'}}\|_1
 \Big\} \quad \mbox{for} \quad 2\le p\le \infty \quad \mbox{and}
 \quad t > 0 \pl .\]
Then there are positive elements $x_1, x_2, x_3, x_4$ with
$x=\sum_k i^k x_k$ and $\|x_k\|_{p,t,d} \le \|x\|_{p,t,d}$.
\end{lemma}

\begin{proof}
Since $\|x^*\|_{p,t,d}=\|x\|_{p,t,d}$, we may clearly assume $x$
self-adjoint. For a self-adjoint element $x$, let $x_+=e_+x$ and
$x_-=e_{-}x$ denote its positive and negative parts, where $e_+$
and $e_-$ stand for the corresponding spectral projections which
commute with $x$. We recall from \cite{Terp} that $L_q(N)$ is a
contractive $N$-bimodule for all $0<q<\infty$. Since $e_+$
commutes with $x$, we obtain
 \[ \|e_+ x\|_{p,t,d} \kl \| x\|_{p,t,d} \pl .\]
The same argument works for $x_- = e_- x$ and the assertion
follows. \qd

At the beginning of section 2 we defined the notion of a
$(q,+)$-summing linear map $T: L_p(N) \to X$. Let $\pi_{q,+}(T)$
denote the infimum of all constants $c$ for which
\eqref{q+summing} holds. The following observation follows
Pisier's argument in \cite{P0}.

\begin{lemma} \label{scha}
Let $N$ be a von Neumann algebra. Let $T: L_p(N) \to X$ be a
$(q,+)$-summing map with $(q,+)$-summing constant $\pi_{q,+}(T)$.
Then, there exists a sequence $(a_n)$ of positive elements of norm
$1$ in $L_p(N)$ such that
 \[ \limm_{n,\U} \big( 1 + \|\frac{T (x_n)}{\pi_{q,+}(T)}\|_X^q \big)^{\frac1q} \kl
  \limm_{n,\U} \|a_n +  x_n\|_p \]
holds for every bounded sequence $(x_n)$ in $L_p(N)_+$ and every
free ultrafilter $\U$.
\end{lemma}

\begin{proof}
Let $C_n$ be the smallest constant satisfying
 \[ \Big( \sum_{k=1}^n \|T (x_k)\|_X^q \Big)^{\frac1q}
 \kl C_n \Big\| \sum_   {k=1}^n x_k \Big\|_p \]
for all families $(x_1, x_2, \ldots, x_n)$ in $L_p(N)_+$. In
particular, we have $\pi_{q,+}(T) = \lim_n C_n$. Let $(\delta_n)$
be a sequence converging to $0$. Then we may find positive
elements $y_1, y_2, \ldots,y_n$ in $L_p(N)$ such that
 \[ \Big( \sum_{k=1}^n \|T(y_k)\|_X^q \Big)^{\frac1q} \lel
 1 \quad \mbox{and} \quad \Big\| \sum_{k=1}^n y_k \Big\|_p \le
 C_n^{-1}(1+\delta_n). \] Let $x_n$ be a
positive element and set $y_{n+1}=x_n$, so that
 \begin{align*}
 \big( 1 + \|T(x_n)\|_X^q \big)^{\frac1q} & =  \Big(
 \sum_{k=1}^{n+1} \|T(y_k)\|_X^q \Big)^{\frac1q} \kl
 C_{n+1} \Big\| \big( \sum_{k=1}^n y_k \big) + x_n \Big\|_p \\
  & \le  \frac{C_{n+1}}{C_n}(1+\delta_n) \, \Big\|
 C_n(1+\delta_n)^{-1} \big( \sum_{k=1}^n y_k \big) + \pi_{q,+}(T)
 x_n \Big\|_p \pl.
 \end{align*}
Let us note that $\|C_n(1+\delta_n)^{-1} \sum_{k=1}^n y_k\|_p \le
1$. Therefore, if we take
 \[ a_n \lel \frac{C_n(1+\delta_n)^{-1} \summ_k y_k}{\ \|C_n(1+\delta_n)^{-1} \summ_k y_k\|_p} \pl ,\]
we obtain
 \[ \big( 1 + \|T(x_n)\|_X^q \big)^{\frac1q} \kl (1+\delta_n)
 \frac{C_{n+1}}{C_n} \, \big\| a_n + \pi_{q,+}(T) x_n
 \big\|_p \quad \mbox{for all} \quad n \ge 1 \pl. \]
Taking the limit yields the assertion.  \qd

\begin{prop}\label{kineq}
Let us fix $2 \le q < p < \infty$. Given a von Neumann algebra $N$
and a $(q,+)$-summing map $T: L_p(N) \to X$, there exists a
sequence of densities $d_n\in L_1(N)$ with $tr(d_n)=1$ such that
 \[ \limm_{n,\U} \|T(x_n)\|_X \kl c(p,q) \, \pi_{q,+}(T)
 \, t^{-\frac{1}{q}} \limm_{n,\U} \|x_n\|_{p,t,d_n} \pl .\]
holds for all $t > 0$ and bounded sequences $(x_n)$ in $L_p(N)$.
In particular, we deduce
\[ \limm_{n,\U} \|T(x_n)\|_X \kl c(p,q) \, \pi_{q,+}(T) \,
\limm_{n,\U} \Big( \|x_n\|_p^{\frac{p'}{q'}} \max \Big\{
\|d_n^{\frac{1}{p'}} x_n\|_1 , \|x_n d_n^{\frac{1}{p'}}\|_1
\Big\}^{1-\frac{p'}{q'}} \Big) \pl .\]
\end{prop}

\begin{proof}
According to Lemma \ref{pos} and the linearity of $T$, we may
clearly assume that the sequence $(x_n)$ lives in the positive cone
$L_p(N)_+$. In particular, according to Lemma \ref{scha} we can find
a sequence $(a_n)$ of norm $1$ positive elements in $L_p(N)$
satisfying
 \[  \limm_{n,\U} \Big( 1 + \big\|
 \frac{T(x_n)}{\pi_{q,+}(T)} \big\|_X^q \Big)^{p/q} \kl
 \limm_{n,\U} \|a_n + x_n\|_p^p \pl .\]
Since $1 + \al \lambda \le
 (1+\lambda)^\al$ for $\lambda > 0$ and $\al > 1 \pl$, we deduce
 \[ 1 + \frac{p}{q} \, \limm_{n,\U} \big\| \frac{T(x_n)}{\pi_{q,+}(T)}
 \big\|_X^q \kl \limm_{n,\U} \|a_n + x_n\|_p^p \pl .\]
Recalling that $a_n$ is norm $1$ in $L_p(N)$,  we obtain by
Theorem \ref{diff}
 \begin{samepage} \begin{eqnarray} \label{mm}
 \limm_{n,\U} \|T(x_n)\|_X
  & \le & \limm_{n,\U} \pi_{q,+}(T) \,
 \Big( \frac{q}{p} \, \big( \|a_n + x_n\|_p^p - \|a_n\|_p^p \big)
 \Big)^{\frac1q} \\ \nonumber
  & \le & \limm_{n,\U} \pi_{q,+}(T) \,
  \Big( q \, 2^{p-1} \, \max \Big\{ \|a_n^{p-1} x_n\|_1 \, , \,
  \|x_n\|_p^p \Big\} \Big)^{\frac1q} \pl .
 \end{eqnarray}\end{samepage}
We define $d_n=a_n^{p}$ and assume that $\|x_n\|_{p,t,d_n}\le 1$.
This implies that
 \[ \max \Big\{ \|a_n^{p-1} x_n\|_1 \, , \, \|x_n\|_p^p \Big\}
 \lel \frac{1}{t} \, \max \Big\{ \big( t^{\frac1p} \|x_n\|_p \big)^p,
 t \|d_n^{\frac{1}{p'}} x_n\|_1, t \|x_n d_n^{\frac{1}{p'}}\|_1
 \Big\} \kl \frac{1}{t} \pl.\]
In conjunction with \eqref{mm}, this proves the first assertion
for sequences $(x_n)$ of positive operators. A further constant
$4$ comes from Lemma \ref{pos} in the general case. Let us prove
the second assertion. We define $\al=\lim_{n,\U}\|x_n\|_p$ and
$\beta=\lim_{n,\U} \max\{\|d_n^{1/p'} x_n\|_1, \|x_n
d_n^{1/p'}\|_1\}$. By the first part we have
 \[ \limm_{n,\U} \|T(x_n)\|_X \kl 4 c(p,q)\pi_{q,+}(T) \inf_{t>0}
 \pl \max \Big\{ t^{\frac1p-\frac1q}\al,t^{1-\frac1q}\beta \Big\}
 \pl .\]
The optimal choice is $t=(\frac{\al}{\beta})^{p'}$ and the optimal
value is then given by
 \begin{align*} \max \Big\{ t^{\frac1p-\frac1q}\al,
 t^{1-\frac1q}\beta \Big\} \lel \al^{\frac{1-1/q}{1-1/p}}
 \beta^{\frac{1/q-1/p}{1-1/p}} \lel
 \al^{p'/q'}\beta^{1-p'/q'} \pl . & \qedhere
 \end{align*} \qd

Our next step through our Nikishin-type result requires some
additional work and in particular the theory of ultraproducts, see
\cite{Ra,RayXu} for some background. Let us assume that $N$ is a
$\si$-finite von Neumann algebra and $d_0$ is a density of a
normal faithful state $\phi_0$. We recall from \cite{Ra} that
 \[ \prodd_{\U}L_p(N) \lel L_p \big( (\prodd_\U N_*)^* \big) \pl .\]
In the following we shall use the notation $(a_n)^{\bullet}$ for
the canonical image of $(a_n)$ in the algebra $(\prod_\U N_*)^*$.
Note that $\big\{ (a_n)^{\bullet} \, | \, \sup_n \|a_n\| < \infty
\big\}$ is dense in $ (\prod_\U N_*)^*$ with respect to the strong
operator topology. Following \cite{RayXu}, we use the support $e$
of the ultraproduct state
 \[ \phi_\U ((a_n)^{\bullet}) \lel \limm_{n,\U} \, tr(d_0 a_n)
 \pl .\]
Let us use the notation $N_\U = e(\prodd_\U N_*)^*e$. Clearly, the
state $\phi_\U$ is a normal faithful state on $N_\U$ and the space
$L_p(N_\U)$ is canonically isomorphic to $e(\prod_{\U} L_p(N))e$,
see \cite{Ra} for further details. This means we can represent
elements $x$ in $L_p(N_\U)$ by sequences of the form
$e(x_n)^{\bullet}e$. This applies in particular for $p=1$ and the
representing sequence for $\phi_\U$ is given by the constant
sequence $(d_0)^{\bullet}$. Here and in the following we also use
the notation $(x_n)^{\bullet}$ for the equivalence class in
$\prod_{\U}L_p(N)$ of a bounded sequence $(x_n)$. Let us recall an
observation from \cite{RayXu}. If $x\in L_p(N)$, then
 \[ (1-e)(x)^{\bullet}\lel 0 \pl .\]
Indeed, we may approximate $x$ by $d_0^{\frac1p}a_n$ with $a_n \in
N$, so that \[ (1-e)(x)^\bullet = (1-e)(d_0^{\frac1p} a_n)^\bullet
\lel 0 \] because $e$ is the support of $(d_0)^{\bullet}$ and
hence the support of its $p$-th root, see again Raynaud's paper
\cite{Ra} for more details on the Mazur map. Our aim is to replace
the sequence of densities $(d_n)$ obtained in Proposition
\ref{kineq} by a single density $d$.

\begin{prop} \label{step3}
Let $2 \le q < p < \infty$ and $\theta=1-\frac{p'}{q'}$. Given a
$\si$-finite von Neumann algebra $N$ and a $(q,+)$-summing map $T:
L_p(N) \to X$, there exists a density $\delta \in L_1(N_\U)$ of a
normal faithful state $\phi_\U$ on $N_\U$ such that
 \[ \limm_{n,\U} \|T(x_n)\|_X \kl c(p,q) \,
 \pi_{q,+}(T) \,
 \|(x_n)^{\bullet}\|_{[L_p(N_\U),
 \Delta_{p,1}(\phi_\U)]_{\theta,1}} \pl. \]
Moreover, if $1 \le r < q$ and $\eta=(\frac{1}{q}-\frac{1}{p}) \big/
(\frac{1}{r}-\frac{1}{p})$, we also have
 \[ \limm_{n,\U} \|T(x_n)\|_X \kl c(p,q) \,
 \pi_{q,+}(T) \,
 \|(x_n)^{\bullet}\|_{[L_p(N_\U),
 \Delta_{p,r}(\phi_\U)]_{\eta,1}} \pl. \]
\end{prop}

\begin{proof}
For the first assertion, it suffices to show that
 \[ \limm_{n,\U} \|T(x_n)\|_X \kl c(p,q) \,
 \pi_{q,+}(T) \, \|(x_n)^{\bullet}\|_p^{1-\theta}
 \|(x_n)^{\bullet}\|_{\Delta_{p,1}(\phi_\U)}^{\theta} \]
for a suitable density $\delta$ of a normal faithful state
$\phi_\U$ in $N_\U$ and $(x_n)^\bullet$ in $L_p(N_\U)$. Indeed,
this is a well-known property of the interpolation bracket $[\, ,
\, ]_{\theta,1}$, see e.g. \cite[p.49]{BL}. On the other hand,
according to Proposition \ref{kineq} we have
$$\limm_{n,\U} \|T(x_n)\|_X \kl c(p,q) \, \pi_{q,+}(T) \,
\|(x_n)^\bullet\|_p^{1-\theta} \limm_{n,\U} \max \Big\{
\|d_n^{\frac{1}{p'}} x_n\|_1, \|x_n d_n^{\frac{1}{p'}}\|_1
\Big\}^{\theta}.$$ Therefore, it remains to find a normal faithful
state $\phi_\U$ for which
 \begin{equation} \label{tarea}
 \limm_{n,\U} \max \Big\{ \|d_n^{\frac{1}{p'}} x_n\|_1,
 \|x_n d_n^{\frac{1}{p'}}\|_1 \Big\} \kl c(p) \,
 \|(x_n)^{\bullet}\|_{\Delta_{p,1}(\phi_\U)}
\end{equation}
whenever $(x_n)^\bullet$ belongs to $L_p(N_\U)$. We deduce for
$(x_n) \in L_p(N_\U)$ that
 \[ \limm_{n,\U} \|d_n^{\frac{1}{p'}}x_n\|_1
 \lel \big\|(d_n^{\frac{1}{p'}})^{\bullet} (x_n)^{\bullet} \big\|_1
 \lel \big\|(d_n^{\frac{1}{p'}})^{\bullet} e (x_n)^{\bullet} \big\|_1
 \lel \big\| \p |(d_n^{\frac{1}{p'}})^{\bullet} e|\pl (x_n)^{\bullet} \big\|_1 \pl
 .\]
Here we use the partial isometry between
$(d_n^{\frac{1}{p'}})^{\bullet} e$ and
$|(d_n^{\frac{1}{p'}})^{\bullet} e|$. Now, we define $$\delta_0
\lel \big( |(d_n^{\frac{1}{p'}})^{\bullet} e|^{2p'} +
(d_0^2)^{\bullet} \big)^{\frac12} \pl.$$ Note that
$$\|\delta_0\|_1 \lel \|\delta_0^2\|_{\frac12}^{\frac12}
\kl \big\| \p |(d_n^{\frac{1}{p'}})^{\bullet} e|^{2p'}
\big\|_{\frac12}^{\frac12} + \big\| (d_0^2)^\bullet
\big\|_{\frac12}^{\frac12} \lel \big\| \p
|(d_n^{\frac{1}{p'}})^{\bullet} e| \p \big\|_{p'}^{p'} + 1 \kl 2
\pl .$$ Thus, if we set $\delta = \delta_0 / \|\delta_0\|_1$, we
obtain the density of a normal faithful state on $N_\U$ given by
$\phi_\U (\cdot) = tr(\delta \, \cdot)$. Indeed, the normality is
clear while the faithfulness follows from the fact that $\delta
\ge \frac{1}{2}(d_0)^{\bullet}$, so that $\delta$ has full
support. It is a state because $\|\delta\|_1 = 1$. We have
 \[ |(d_n^{\frac{1}{p'}})^{\bullet} e|^{2p'} \le \delta_0^2 \
 \Rightarrow |(d_n^{\frac{1}{p'}})^{\bullet} e| \le
 \delta_0^{\frac{1}{p'}} \pl.  \]
Hence we can find a contraction $w$ in $N_\U$ such that
$|(d_n^{\frac{1}{p'}})^{\bullet} e| = w \delta_0^{\frac{1}{p'}} =
\delta_0^{\frac{1}{p'}} w$. This implies
$$\big\| \p |(d_n^{\frac{1}{p'}})^{\bullet} e| \p (x_n)^{\bullet}
\big\|_1 \lel \big\| w \delta_0^{\frac{1}{p'}} (x_n)^{\bullet}
\big\|_1 \kl \big\| \delta_0^{\frac{1}{p'}}(x_n)^{\bullet}
\big\|_1 \kl 2^{\frac{1}{p'}} \big\|
\delta^{\frac{1}{p'}}(x_n)^{\bullet} \big\|_1 \pl.$$ Similarly, we
have $$\big\| (x_n)^{\bullet} (d_n^{\frac{1}{p'}}) \big\|_1 =
\big\| (x_n)^{\bullet} e (d_n^{\frac{1}{p'}}) \big\|_1 = \big\|
(x_n)^{\bullet} \p |(d_n^{\frac{1}{p'}})^\bullet e| \p \big\|_1 =
\big\| (x_n)^{\bullet} \delta_0^{\frac{1}{p'}} w \big\|_1 \le
2^{\frac{1}{p'}} \big\| (x_n)^{\bullet} \delta^{\frac{1}{p'}}
\big\|_1$$ for all $(x_n)\in L_p(N_\U)$. Therefore, we obtain
\eqref{tarea} and the first assertion is proved. The second
assertion is an immediate consequence of the first one and the
reiteration theorem. Indeed, according to Theorem \ref{minint} we
have
 \[ \Delta_{p,r}(\phi_\U) \lel \big[ L_p(N_\U),
 \Delta_{p,1}(\phi_\U) \big]_{\zeta} \]
where $1/r = \zeta + (1-\zeta)/p$, so that $\zeta =
1-\frac{p'}{r'}$. The reiteration theorem for the real method
\cite[Theorem 4.7.2]{BL} implies $[L_p(N_\U),
\Delta_{p,r}(\phi_\U)]_{\eta,1} = [L_p(N_\U),
\Delta_{p,1}(\phi_\U)]_{\theta,1}$ with $\theta=\eta \zeta$. We find
$\eta=(\frac{1}{q}-\frac{1}{p}) \big/ (\frac{1}{r}-\frac{1}{p})$ as
announced. \qd

\begin{cor} \label{step5}
Let $2 < q < p < \infty$ and $\eta = (\frac1q - \frac1p) \big/
(\frac12 - \frac1p)$. Given any von Neumann algebra $N$ and a
$(q,+)$-summing map $T: L_p(N) \to X$, there exists a density $d \in
L_1(N)$ with $tr(d) = 1$ and support $e$ such that
 \[ \|T(x)\|_X \kl c(p,q) \, \pi_{q,+}(T) \,
 \|x\|_{[eL_p(N)+L_p(N)e, \Delta_{p,2}(\phi)]_{\eta,1}} \pl .\]
\end{cor}

\begin{proof}
Let us first assume that $N$ is a $\si$-finite von Neumann algebra
and set $\frac{1}{r}=\frac{1}{2}-\frac{1}{p}$. We use the density
$\delta = (\delta_n)^\bullet \in L_1(N_\U)$ from Proposition
\ref{step3}~. Given $x \in L_p(N)$ we observe that we have
 \[ \limm_{n,\U} \|\delta_n^{\frac{1}{r}} x\|_2^2
 \lel \limm_{n,\U} tr \big( \delta_n^{\frac{1}{r}} xx^*
 \delta_n^{\frac{1}{r}} \big) \lel \limm_{n,\U} tr \big(
 xx^* \delta_n^{\frac{2}{r}} \big) \pl .\]
This defines a positive element $a = \lim_{n,\U}
\delta_n^{\frac2r} \in L_{\frac{r}{2}}(N)$. We take
$d=a^{\frac{r}{2}}$ and recall that
 \[ \limm_{n,\U} \|x\|_p^{1-\eta} \max \Big\{
 \|\delta_n^{\frac{1}{r}} x\|_2, \|x
 \delta_n^{\frac{1}{r}}\|_2 \Big\}^{\eta} \lel \|x\|_p^{1-\eta}
 \max \Big\{ \|d^{\frac{1}{r}} x\|_2, \|x
 d^{\frac{1}{r}}\|_2 \Big\}^{\eta} \pl .\]
Thus Proposition \ref{step3} applied to the constant sequence
$(x)^{\bullet}$ yields the result, because
$(x)^{\bullet}=e(x)^{\bullet}e$. When $N$ is an arbitrary von
Neumann algebra, we choose a normal strictly semifinite weight
$\psi=\lim_i \phi_i$ such that $\phi_i$ is a positive functional
and the support $e_i$ of $\phi_i$ satisfies
$\si_t^{\phi_i}(e_i)=e_i$. Then $N_i=e_iNe_i$ is $\si$-finite and
we find a density $d_i\in L_1(N_i)$ with $tr(d_i)=1$ and such that
 \[ \|T(e_ixe_i)\|_X \kl c(p,q) \pl \pi_{q,+}(T) \,
 \|e_ixe_i\|_p^{1-\eta} \max \Big\{ \|d_i^{\frac{1}{r}} x\|_2,
 \|x d_i^{\frac{1}{r}}\|_2 \Big\}^{\eta} \pl .\]
As above we can pass to the limit $d^{\frac{2}{r}} = \lim_i
d_i^{\frac{2}{r}}$.  \qd

\begin{rem}{\rm It is tempting to use a weak limit $d=\lim_{n,\U}
d_n^{1/p'}\in L_{p'}(N)$ in Proposition \ref{kineq}. The problem we
face is the equality
 \begin{equation}\label{limitq}
  \|dx\|_1 \pl \stackrel{?}{=} \pl  \limm_{n,\U} \|d_n^{1/p'}x\|_1
  \pl .
  \end{equation}
This equality does not hold in general. Indeed, assuming
\eqref{limitq}, we would deduce from the polar decomposition that
$\|w^*-\lim_{n,\U} a_n x\|_1 = \pl \lim_{n,\U} \|a_nx\|_1$ holds
for all bounded sequences $a_n$ and $x\in
S_{p}=L_p(B(\ell_2),tr)$. In $S_{p'}$ we may choose $a_n=e_{n,1}$.
Then $\lim_n a_n=0$ weakly and $\lim_n \|a_ne_{11}\|_1 \lel 1$. We
suspect that we need some equi-integrability for \eqref{limitq} to
hold. Our proof does not provide any equi-integrability condition.
}\end{rem}

\noindent We are ready for the main result.

\begin{theorem} \label{Nik-th}
Let $2\le q <p<\infty$ and let $N$ be any von Neumann algebra.
Given a $(q,+)$-summing map $T: L_p(N) \to X$, there exists a
density $d \in L_1(N)$ such that the inequality below holds for
any index $q<r<p$
 \[ \|T(x)\|_X \kl c(p,q,r) \pl \pi_{q,+}(T) \pl
 \|x\|_{\Delta_{p,r}(\phi)} \pl .\]
Moreover, there exists a map $\widetilde{T}: L_r(N \oplus N) \to
X$ such that
 \[ T(x) \lel \widetilde{T} \big( d^{\frac1r-\frac{1}{p}} x,
 x d^{\frac1r-\frac{1}{p}} x \big) \lel \widetilde{T} \big(
 j_{p,r}(x) \big) \quad \mbox{and} \quad \|\widetilde{T}\|
 \kl c(p,q,r) \, \pi_{q,+}(T) \pl. \]
\end{theorem}

\begin{proof}
According to Corollary \ref{step5}, we may find $d$ such that
 \[ \|T(x)\|_X \kl c(p,q) \pl \pi_{q,+}(T) \, \|x\|_p^{1-\eta}
 \max \Big\{ \|d^{\frac12-\frac1p} x\|_2, \|x
 d^{\frac12-\frac1p}\|_2 \Big\}^{\eta} \pl .\]
Let $e$ denote the support of $d$ and let $\phi(x) = tr(dx)$ be
the associated state. Let us decompose any element $x$ in $L_p(N)$
as $x = exe+ex(1-e)+(1-e)xe+(1-e)x(1-e)$. Then we note that $T$
vanishes on the corner $(1-e)L_p(N)(1-e)$. We apply Lemma
\ref{fact} for $\tilde{\eta}<\eta$ and deduce
 \begin{equation}\label{pp} \|T(x)\|_X \le c(p,q) c(\eta,\tilde{\eta}) \pl
 \pi_{q,+}(T) \, \|x\|_{[eL_p(N)+L_p(N)e,\Delta_{p,2}(\phi)]_{\tilde{\eta}}}
 \end{equation}
for all $x\in L_p(N)$ such that $(1-e)x(1-e)=0$. We recall the
isomorphism from  Lemma \ref{three}:
  \begin{equation}\label{pppp}
 [eL_p(N)+L_p(N)e,\Delta_{p,2}(\phi)]_{\tilde{\eta}}
 \simeq  e\Delta_{p,r}(\phi)e \oplus eL_r(N)(1-e) \oplus
 (1-e)L_r(N)e \pl ,
 \end{equation}
where $\frac1r = \frac{1-\tilde{\eta}}{p} +
\frac{\tilde{\eta}}{2}=\frac1p+\tilde{\eta}(\frac12-\frac1p)<
\frac1p+\eta(\frac12-\frac1p)=\frac1q$. Thus for every $q<r<p$ we
can find a suitable $\tilde{\eta}$. Note that $\eta = 1$ when $q=2$,
but then we may use that $(2,+)$-summing implies $(q,+)$-summing for
all $q > 2$. We denote by
 \[ T_{\tilde{\eta}}: e\Delta_{p,r}(\phi)e \oplus eL_r(N)(1-e) \oplus
 (1-e)L_r(N)e \to X\]
the corresponding bounded map given by \eqref{pp} and
\eqref{pppp}. Now we proceed as in Theorem \ref{p} and define
 \[ \widetilde{\mathcal{Q}}_r:L_r(N)\oplus L_r(N)\to
 e\Delta_{p,r}(\phi)e \oplus eL_r(N)(1-e) \oplus
 (1-e)L_r(N)e\]
by $\widetilde{\mathcal{Q}}_r(x,y)\lel \mathcal{Q}_r(exe,eye)\oplus
ex(1-e)\oplus (1-e)ye$, where $\mathcal{Q}_r$ is the projection from
Theorem \ref{minint}.  This allows us to define
$\widetilde{T}(x,y)\lel T_{\tilde{\eta}}\widetilde{\mathcal{Q}}_r$
so that
 \begin{align*}
  \widetilde{T}(d^{\frac1r-\frac1p}x,xd^{\frac1r-\frac1p})
 &=  T_{\tilde{\eta}} \big((exe,exe) \oplus
 d^{\frac1r-\frac1p}x(1-e)\oplus (1-e)xd^{\frac1r-\frac1p}\big) \\
 &= T(exe+(1-e)x+ex(1-e))\lel T(x) \pl .
 \end{align*}
The norm estimate follows from Theorem \ref{minint}, \eqref{pp} and
\eqref{pppp} (see Lemma \ref{three}).\qd


\begin{rem}{\rm  In the commutative analog of Theorem \ref{Nik-th},
the restriction $p,q \ge 2$ is not needed. It would be very
interesting to know whether this restriction is necessary in the
noncommutative setting.}
\end{rem}

\begin{rem}{\rm We know from Lust-Piquard work \cite{lust-gro} that
for $q=2$ the situation is much nicer. Let $\frac1q + \frac1s =
\frac12$. Then, the noncommutative Khintchine inequality \cite{LP}
implies
 \[ \summ_k \|T(x_k)\|_X^2\kl c_p^2 \pl c_2(T) \pl \sup_{\|a\|_{s/2}\le 1}
 \sum_k tr \big( a(x_k^*x_k+x_kx_k^*) \big) \pl. \]
Applying the standard separation argument one  obtains a
factorization $T=vj_{p,2}$ through the inclusion map $id:L_p(N)\to
\Delta_{p,2}(\phi)$. We refer to \cite{lust-gro} for more details
and to \cite{LPPX} for further information. 
}\end{rem}

\begin{cor}\label{thmB}
Let $1< q<r\le 2$ and $T:X\to L_q(N)$ be a linear map such that
$T^*$ has cotype $r'$. Then there exists a density $d\in L_1(N)$
with $tr(d)=1$ such that for every $q<p<r$ there exists a map
$u:X\to L_p(N)$ satisfying
 \[ T(x) \lel d^{\frac1q-\frac1p} u(x) + u(x) d^{\frac1q-\frac1p}
 \quad \mbox{for all} \quad x\in X \pl . \]
\end{cor}

\begin{proof} By Theorem \ref{Nik-th} we find a density
$d$ with $tr(d)=1$ such that
 \[ \big\| T^*: \Delta_{q',p'}(\phi)\to X^* \big\| \kl
 c(p',q',r') \pl \pi_{r',+}(T^*) \kl  2 \pl c(p',q',r') \pl
 c_{r'}(T^*) \]
for every $r'<p'<q'$. We write
$T^*_{q',p'}:\Delta_{q',p'}(\phi)\to X^*$ for the corresponding
map. Let $\widetilde{\mathcal{Q}}_{p'}: L_{p'}(N) \to
\Delta_{q',p'}(\phi)$ be the projection from Theorem \ref{Nik-th}.
We recall that $\mathcal{Q}_{p'}^*=(v_e,v_e)$ has two identical
components. As in the proof of Theorem \ref{p} we define
$v:\Delta_{q',p'}(\phi)^*\to L_p(N)$ by
$v(\xi_e,\xi_{1-e,e},\xi_{e,1-e})=v_e(\xi_e)+\xi_{1-e,e}+\xi_{e,1-e}$.
Following the argument from  Theorem \ref{p} we can check that $u
=v (T^*_{q',p'})^*: X \to L_{p}(N)$ provides the corresponding
decomposition. Note that in the $\si$-finite case we may assume
that $d$ has full support. Then formally
 \[ u(x) \lel (\L_{d^{\frac1q-\frac1p}}
 + \R_{d^{\frac1q-\frac1p}})^{-1}T(x) \pl .\]
In full generality we have
 \begin{align*}
  u(x) &=
 (\L_{d^{\frac1q-\frac1p}} + \R_{d^{\frac1q-\frac1p}})^{-1}
 eT(x)e + d^{\frac1p-\frac1q}T(x)(1-e)+ (1-e)T(x)d^{\frac1p-\frac1q}
 \pl .
 \end{align*}
We should warn the reader that these multiplications are usually
not well-defined, see both Step 2 and Step 3 of the proof of
Theorem \ref{minint} for a rigorous interpretation using
Haagerup's construction. \qd

\begin{proof}[{\bf Proof of Theorem \ref{mainq}.}]
 Let $1\le q< 2$ and let $X$ be an infinite-dimensional
 subspace of $L_q(N)$ not containing $\ell_q$.
 According to Raynaud and Xu's result  \cite[Theorem5.1]{RayXu} we
 deduce that $X$ does not contain $\ell_q(n)$'s uniformly. By
 the  Krivine-Maurey-Pisier  theorem, the type index of $X$ satisfies
 $p_X>q$. Let $q < r < p_X$ so that $X^*$ has cotype $r'$. Let
 $\iota:X\to L_q(N)$ be the inclusion map. Then
 $T=\iota^*:L_{q'}(N)\to X^*$ has cotype $r'$ and the
 assertion follows from Corollary \ref{thmB}.
\qd

\noindent We refer to \cite{Ran,RayXu} for the definition of
$q$-equiintegrable sets in $L_q(N)$.

\begin{cor}
If $X \subset L_q(N)$ and $1 \le q < 2$, the following are
equivalent
\begin{enumerate}
\item[i)] The unit ball of $X$ is $q$-equiintegrable.

\item[ii)] There exists a density $d\in L_1(N)$ such that
$$u: x \in X \mapsto \big( d^{\frac1r-\frac1q} x,
x d^{\frac1r-\frac1q} \big) \in L_r(N\oplus N)$$ is an isomorphic
embedding for some $(all)$ $0<r<q$.

\item[iii)] There exists $q<p<2$ and a bounded linear map $$u: X
\to L_p(N)$$ such that $x = d^{\frac1q-\frac1p} u(x) + u(x)
d^{\frac1q-\frac1p}$ for some positive density $d \in L_1(N)$.
\end{enumerate}
\end{cor}

\begin{proof} According to \cite[Theorem 5.1]{RayXu}, the conditions
i) and ii) are both equivalent to the fact that $X$ does not
contain $\ell_q$ and hence imply iii) by means of Theorem
\ref{mainq}. On the other hand, iii) implies that $X$ has type
$p>q$ and hence can not contain $\ell_q$.\qd

\noindent \textbf{Acknowledgement.} The referee's comments of the
first version of this paper not only improved the presentation,
but enabled us to simplify the results considerably. We are also
indebted to him for allowing  us to present his proof of Theorem
\ref{pr}. We would like to thank Narcisse Randrianantoanina for
keeping us up to date on his closely related work. The first-named
author was partially supported by the National Science Foundation,
DMS-0301116 and DMS 05-56120. The second-named author was
partially supported by \lq\lq Programa Ram{\'o}n y Cajal, 2005" and
the Grant MTM2004-00678, Spain.

\vskip0.3cm

\noindent  {\bf  2000 Mathematics Subject Classification:} Primary
46L53, Secondary 46B25.

\noindent {\bf Keyword's:} Rosenthal's theorem, noncommutative
$L_p$ spaces.

\vskip0.3cm

\noindent \textbf{Marius Junge} \\
\textsc{Department of Mathematics} \\
\textsc{University of Illinois at Urbana-Champaign} \\ 273 Altgeld
Hall, 1409 W. Green Street, Urbana, IL 61801, USA \\
\texttt{junge@math.uiuc.edu}

\vskip0.3cm

\noindent \textbf{Javier Parcet} \\
\textsc{Departamento de Matem{\'a}ticas} \\
\textsc{Instituto de Matem}{\scriptsize {\'A}}\textsc{ticas
y F}{\scriptsize {\'I}}\textsc{sica Fundamental} \\
\textsc{Consejo Superior de Investigaciones Cient}{\scriptsize
{\'I}}\textsc{ficas} \\ Depto de
Matem{\'a}ticas, Univ. Aut{\'o}noma de Madrid, 28049, Spain \\
\texttt{javier.parcet@uam.es}

\end{document}